\documentclass[12pt]{article}
\usepackage{amsmath}
\usepackage{amssymb}
\usepackage{rotating}
\usepackage{latexsym}
\usepackage[all]{xy}
\newenvironment{conjalt}[1]{\vspace{1ex}
\par\noindent{\bf Conjecture} #1 \em}{\vspace{1ex}\par}
\newtheorem{lemma}{Lemma}
\newtheorem{prop}{Proposition}

\newcounter{example}
\newcounter{rem}%counter for \rem (Remark) command

\newcommand{\barbbQ}{\ensuremath{\bar{\bbQ}}}

\newcommand{\bbC}{\ensuremath{\mathbb{C}}}

\newcommand{\bbQ}{\ensuremath{\mathbb{Q}}}

\newcommand{\bbR}{\ensuremath{\mathbb{R}}}
\newcommand{\bbZ}{\ensuremath{\mathbb{Z}}}

\newcommand{\beq}{\begin{equation}}
\newcommand{\beql}[1]{\begin{equation}\label{#1}}
\newcommand{\bPf}{\noindent \textsc{Proof\ }}

\newcommand{\cH}{\ensuremath{\mathcal{H}}}

\newcommand{\cO}{\ensuremath{\mathcal{O}}}

\newcommand{\cR}{\ensuremath{\mathcal{R}}}

\newcommand{\cU}{\ensuremath{\mathcal{U}}}

\newcommand{\displaymapdef}[5]
{\[
\begin{array}{rcrcl}
 #1 &:& #2 &\longrightarrow& #3 \\
    & &    &                    \\
    & & #4 &\longmapsto    & #5
\end{array}
\]}

\newcommand{\eeq}{\end{equation}}
\newcommand{\eg}{\emph{e.g.}}
\newcommand{\eitheror}[4] %% to be used to create RHS of an equation, ineq. etc. with 2 alternatives and their conditions
{\left\{                  %% arguments #1 and #3 are the alternative RHS's (in mathmode) and
\begin{array}{ll}         %% arguments #2 and #4 are the respective conditions (in paragraph mode)
#1& \mbox{#2}\\
#3& \mbox{#4}
\end{array}
\right.}
\newcommand{\ePf}{\hspace*{\fill}~$\Box$\vertsp\par}%ends proof
\newcommand{\eps}{\ensuremath{\varepsilon}}

\newcommand{\ff}{\ensuremath{\mathfrak{f}}}

\newcommand{\fp}{\ensuremath{\mathfrak{p}}}
\newcommand{\fP}{\ensuremath{\mathfrak{P}}}
\newcommand{\fq}{\ensuremath{\mathfrak{q}}}
\newcommand{\fQ}{\ensuremath{\mathfrak{Q}}}

\newcommand{\fs}{\ensuremath{\mathfrak{s}}}
\newcommand{\fS}{\ensuremath{\mathfrak{S}}}

\newcommand{\Gal}{{\rm Gal}}

\newcommand{\half}{\frac{1}{2}}

\newcommand{\ie}{\emph{i.e.}}
\newcommand{\im}{\ensuremath{{\rm im}}}
\newcommand{\Ind}{{\rm Ind}}
\newcommand{\inv}{^{-1}}% use only in mathmode

\newcommand{\ndiv}{\nmid}
\newcommand{\nin}{\not\in}
\newcommand{\ord}{{\rm ord}}

\newcommand{\rec}{{\rm rec}}

\newcommand{\rem}{\refstepcounter{rem}\noindent{\sc Remark \therem}}%this
\newcommand{\ses}[5]{#1\rightarrow #2\longrightarrow #3\longrightarrow #4\rightarrow #5}

\newcommand{\vertsp}{\vspace{1ex}}

\newcommand{\pnpo}{\ensuremath{{p^{n+1}}}}

\newcommand{\sKk}{\ensuremath{\fs_{K/k}}}
\newcommand{\SKk}{\ensuremath{\fS_{K/k}}}
\newcommand{\SKkS}{\ensuremath{\fS_{K/k,S}}}
\newcommand{\sKkS}{\ensuremath{\fs_{K/k,S}}}
\newcommand{\Sram}{\ensuremath{S_{\rm ram}}}

\newcommand{\WUoKp}{\ensuremath{{\textstyle \bigwedge^d_{\bbZ_pG}U^1(K_p)}}}

\newcommand{\WUSKp}{\ensuremath{{\textstyle \bigwedge^d_{\bbQ \bar{G}}\bbQ U_S(K^+)}}}
\newcommand{\WUotwo}{\ensuremath{{\textstyle \bigwedge^2_{\bbZ_pG}U^1(K_p)}}}

\newcommand{\Zbp}{\ensuremath{\bbZ_{(p)}}}
\newcommand{\ZpnpoZ}{\ensuremath{\bbZ/\pnpo \bbZ}}

\newcommand{\ZpnpoZst}{\ensuremath{(\ZpnpoZ)^\times}}

\newcommand{\onecolc}[2]{\multicolumn{#1}{|c|}{#2}}

\newcommand{\phv}[2]{\vrule height #1pt depth #2pt width 0pt}
\begin{document}
\title{Verifying the Congruence Conjecture\\ for Rubin-Stark Elements}
\author{X.-F. Roblot\\
Universit\'e de Lyon, Universit\'e Claude Bernard Lyon 1 \\
Institut Camille Jordan, CNRS -- UMR 5208 \\
roblot@math.univ-lyon1.fr\vertsp\\D. Solomon
\\King's College, London\\
david.solomon@kcl.ac.uk}
\maketitle

\begin{abstract}
  The `Congruence Conjecture' was developed by the second author in a
  previous paper~\cite{laser}. It provides a conjectural explicit
  reciprocity law for a certain element associated to an abelian extension of a totally real number field
  whose existence is predicted by earlier conjectures of Rubin and Stark.
  The first aim of the present paper is to design and apply techniques to investigate the Congruence Conjecture numerically.
  We then present complete verifications of the conjecture in 48 varied cases with real quadratic base fields.

\end{abstract}

\section{Introduction}

The primary purpose of this paper is to provide numerical evidence for the `Congruence Conjecture'. This first appeared
as Conjecture~5.4 of~\cite{twizo} but we shall refer here to the improved and generalised version appearing
as $CC(K/k,S,p,n)$
in~\cite{laser}.
Thus $K/k$ denotes an
abelian extension of number fields, $S$ a finite set of places of $k$, $p$ an {\em odd} prime number and $n$ an integer, $n\geq -1$.
We suppose that $k$ is totally real of degree $d$ and that $K$ is of $CM$ type and contains $\xi_\pnpo:=\exp(2\pi i/\pnpo)$.
(More precise conditions on $S$ will be explained later.)
In this set-up, we can say that the $CC$ is a conjectural,
{\em $p$-adic, explicit reciprocity law
for the so-called Rubin-Stark element $\eta_{K^+/k,S}$}. We recall that $\eta_{K^+/k,S}$ is a particular element of
a certain $d$th exterior power of the global $S$-units of $K^+$ (tensored with $\bbQ$)
which is predicted to exist by Stark's
conjectures, as reformulated and refined by Rubin in~\cite{Rubin}. It is uniquely
determined by the $d$th derivatives at $s=0$ of the
$S$-truncated Artin $L$-functions of even characters of $\Gal(K/k)$.

By way of illustration, consider the simplest case $K/k=\bbQ(\xi_\pnpo)/\bbQ$, $S=\{\infty,p\}$ (so $d=1$). One can then prove that
$\eta_{K^+/k,S}$ exists and equals $-\half\otimes(1-\xi_\pnpo)(1-\xi_\pnpo\inv)$. Moreover,
the $CC$ then reduces to the explicit reciprocity law
proven by Artin and Hasse in~\cite{A-H}. This is
a precise formula for the Hilbert symbol $(1-\xi_\pnpo,u)_{K_\fP,\pnpo}$, for any $u\in U^1(K_\fP)$, which involves the $p$-adic logarithms of
the conjugates of $u$ over $\bbQ_p$.
(Here $K_\fP$ denotes the completion of
$K$ at the unique prime $\fP$ dividing $p$ and $U^1(K_\fP)$ its group of principal units.)

For the general case of the $CC$ one must replace $u$ by an element $\theta$ of a certain $d$th exterior power of $U^1(K_p)$
(the principal, $p$-semilocal units of $K$). From  $\theta$ and $\eta_{K^+/k,S}$ one forms a $d\times d$ determinant
of (additive, group-ring-valued) Hilbert symbols.
The conjectural reciprocity law takes the form of a congruence modulo $\pnpo$  between this determinant and
$\sKkS(\theta)$ (for any $\theta$),
where $\sKkS$ is a map defined explicitly in~\cite{laser} and~\cite{twizo} using a certain $p$-adic regulator and the values at $s=1$ of the
$S$-truncated Artin $L$-functions of odd characters of $\Gal(K/k)$.
More details of $\eta_{K^+/k,S}$, this determinant,
the map $\sKkS$ and the precise formulation of $CC(K/k,S,p,n)$ are given in Section~\ref{sec: the $CC$}.

In the case where $K$ is absolutely abelian, the
$CC$ was proven (with some restrictions) in~\cite{laser}: One reduces first to the case $k=\bbQ$
where $\eta_{K/k,S}$ is essentially a cyclotomic unit (as above) and the $CC$ can be proven without restriction, replacing the Artin-Hasse
law with a generalisation due to Coleman.
This case of the $CC$ (or more precisely the connection it makes
between reciprocity laws and the map $\fs_{K/\bbQ,S}$) finds applications in Iwasawa Theory related to some new annihilators
of the class-groups of real abelian fields (see~\cite{sonic}).
This gives one motivation for
studying the $CC$ more generally.

Unfortunately, there are very few cases with $K$ not absolutely abelian in which
$CC(K/k,S,p,n)$ can be proven, even partially (see~\cite[\S 4]{laser}). Indeed for such $K$,
one can't even prove the {\em existence} of $\eta_{K^+/k,S}$ except in very special cases (see Section~\ref{sec: RS + pairing H}).
On the other hand, techniques for the {\em numerical} computation of $\eta_{K^+/k,S}$ were developed by the authors in~\cite{Solomon-Roblot1}.
A slight simplification of these methods is used in the present paper to identify $\eta_{K^+/k,S}$ with virtual certainty. The rest of the paper
is concerned with the detailed numerical verification of 48 varied cases of the $CC$ using the computed values of $\eta_{K^+/k,S}$.

In order to make the computations manageable we still need to restrict the parameters $(K/k,S,p,n)$: in all our test cases $k$ is (real)
quadratic, $p\leq 7$ and $n=0$ or $1$. (On the other hand, $K/\bbQ$
is always non-abelian and frequently non-Galois). The precise set-up is given
at the start of Section~\ref{sec: methods of comp}. We then explain in detail how we computed
the objects appearing in the $CC$, in order:
the map $\sKkS$, economical sets of (Galois) generators for $U^1(K_p)$ and its exterior square,
the element $\eta_{K^+/k,S}$ and the Hilbert-symbol-determinant $H_{K/k,n}(\eta_{K^+/k,S},\theta)$.
Some of our techniques are well known and even implemented in PARI/GP (which is also the medium of all our computations).
However, we believe that the majority are
innovative and may well find applications elsewhere.
It is worth mentioning an important dichotomy which emerges in our examples, between the minority of cases in which $p$ divides $[K:k]$
and the majority in which it does not.
On the one hand, the former cases provide a more probing test of the
conjecture. For instance, since $k$ is quadratic, the condition $n=1$ requires $p|[K:k]$. On the other hand, cases of the latter type are
much quicker to compute.

Finally, Section~\ref{sec: results} presents the results of the computations.
One simple but characteristic example is explained in detail. Data from the
remaining ones are summarised in tables at the end of the paper.

Some notations and conventions: All number fields are finite extensions of $\bbQ$ within $\barbbQ$ which is
the algebraic closure of $\bbQ$ within $\bbC$. If $F$ is any field and $m$ any positive integer, we shall write $\mu_m(F)$ for the group of all
$mth$ roots of unity in $F$. We shall abbreviate $\mu_m(\bbC)$ to $\mu_m$ and write $\xi_m$ for its generator $\exp(2\pi i/m)$.
Suppose $L/F$ is a Galois extension of number fields and $\fQ$ a prime ideal of $\cO_L$ with $\fq=F\cap\fQ$. We shall write $D_\fQ(L/F)$ for the
decomposition subgroup of $\Gal(L/F)$ at $\fQ$ and similarly $T_\fQ(L/F)$ for the inertia subgroup.
We shall identify $D_\fQ(L/F)$ with the Galois group
of the completed extension $L_\fQ/F_\fq$ and $T_\fQ(L/F)$ with its inertia group in the usual way. If $\cR$ is a commutative ring and $H$ is a
finite group, we shall write simply $\cR H$ for the group-ring often denoted $\cR[H]$.

The second author wishes to thank Cristian Popescu and UCSD for their hospitality during the sabbatical year in
which part of this paper was written.

\section{The Congruence Conjecture}\label{sec: the $CC$}
\subsection{The Map $\sKkS$}\label{subsec: TheMapsKkS}
Given an abelian extension $K/k$ of number fields as above, we
write $G$ for $\Gal(K/k)$ and $S_\infty=S_\infty(k)$  and $\Sram=\Sram(K/k)$
respectively for the  set of infinite places of $k$ and the set of those
finite places of $k$ which ramify in $K$. We always identify finite
places with prime ideals so, for instance, $\Sram$ consists of the
prime factors of the conductor $\ff(K)$ of $K/k$. We denote by $S_p=S_p(k)$ the set of places of $k$ dividing
the prime number $p\neq 2$. The finite $S$ appearing in the $CC$ must satisfy the hypothesis
%%For the purposes of this paper
%% we shall assume the following hypotheses (and notations) concerning $K/k$, $p$, an integer $n\geq 0$ and a finite set $S$ of places of $k$ :
\beq\label{subhyp: S contains...}
\mbox{$S$ contains $S^1:=S_\infty\cup\Sram\cup S_p$}
\eeq
which we assume henceforth. Recall also that we are assuming $K$ is CM so that $[K:K^+]=2$ where $K^+$ is its maximal real subfield
which contains $k$. The extra assumption that $K$ contains $\mu_{\pnpo}$, which is necessary for the $CC$, may be dropped until further notice.

If $s$ is a complex number with ${\rm Re}(s)>1$, we define an Euler product in the complex group ring $\bbC G$ of $G$ by
\beq\label{eq: introduction}
\Theta_{K/k, S}(s):=\prod_{\fq\nin S}\left(1-N\fq^{-s}\sigma_\fq\inv\right)\inv\
\eeq
(The prime ideal $\fq$ of $\cO_k$ ranges over all those not in $S$ and $\sigma_\fq$ denotes the corresponding Frobenius element of $G$.)
Indeed, the condition  ${\rm Re}(s)>1$ implies that
$(1-N\fq^{-s}\sigma_\fq\inv)$ lies in $\bbC G^\times$ and that the product converges absolutely. $\Theta_{K/k, S}(s)$ is sometimes called the
{\em `equivariant $L$-function'} because, if $\hat{G}$ denotes the group of (complex, irreducible) characters of $G$, then one can write
$
\Theta_{K/k, S}(s)
=\sum_{\chi\in \hat{G}}L_{K/k, S}(s,\chi)e_{\chi\inv,G}
$.
Here, for any $\chi\in \hat{G}$, we write
$e_{\chi,G}$ for the corresponding
idempotent $\frac{1}{|G|}\sum_{g\in G}\chi(g)g\inv$ of $\bbC G$ and $L_{K/k, S}(s,\chi)$ for the ($S$-truncated Artin) $L$-function, \ie\ the function whose
Euler product for ${\rm Re}(s)>1$ is obtained by applying $\chi\inv$ termwise to the R.H.S. of~({\ref{eq: introduction}}). Since
$L_{K/k, S}(s,\chi)$ extends to a meromorphic function on $\bbC$ so does $\Theta_{K/k, S}(s)$
(with values in $\bbC G$).
Now let $c$ denote
the element of $G$ determined by complex conjugation, so that $\Gal(K/K^+)=\{1,c\}=\langle c\rangle$. A character $\chi\in\hat{G}$ is called
{\em odd} ({\em resp.} {\em even}) if and only if $\chi(c)=-1$ ({\em resp.}
$\chi(c)=1$). If $\cR$ is any
commutative ring in which $2$ is invertible, we write $e^\pm$ for the two idempotents $\half(1\pm c)\in\cR \langle c\rangle$.
Any $\cR \langle c\rangle$-module $M$ then splits as $M^+\oplus M^-$ where $M^+=e^+ M$ is the `plus-submodule' and $M^-=e^- M$ is the `minus-submodule'.
Taking $\cR=\bbC$ and $M=\bbC G$, we get a corresponding decomposition
$\Theta_{K/k, S}(s)=e^+\Theta_{K/k, S}(s)+e^-\Theta_{K/k, S}(s)=:\Theta^+_{K/k, S}(s)+\Theta^-_{K/k, S}(s)$, say.  Clearly,
$\Theta^-_{K/k, S}(s)=\sum_{\chi\ {\rm odd }}L_{K/k, S}(s,\chi)e_{\chi\inv,G}$  and since $L_{K/k, S}(s,\chi)$ is regular
at $s=1$ whenever $\chi$ is not the trivial character $\chi_0$, it follows that $\Theta^-_{K/k, S}(s)$ is also regular there. We set
\beq\label{eq: defn of a}
a^-_{K/k,S}:=\left(\frac{i}{\pi}\right)^d \Theta^-_{K/k,S}(1)=\left(\frac{i}{\pi}\right)^d\sum_{\chi\in\hat{G}\atop\chi\ {\rm odd}}L_{K/k, S}(1,\chi)e_{\chi\inv,G}
\eeq
In this notation, it is not hard to see that $a^-_{K/k,S}$ lies in $i^d\bbR G^-$. In fact it lies in
$\barbbQ G^-$ and indeed a much finer statement will be proven in Proposition~\ref{prop: stuff on atilde}.
%% follows from~\cite[Prop 3.1]{twizo} (and~\cite[eq. (9)]{laser}),
%% that the coefficients
%% of $a^-_{K/k,S}$ lie in $\barbbQ$, and indeed in $\sqrt{d_k}\bbQ(\mu_{f(K)})$
%% where $d_k$ denotes the absolute discriminant of $k$ and $f(K)$ is the positive generator
%% of the ideal $\ff(K)\cap\bbZ$.

For each $\fP\in S_p(K)$ we write $K_\fP$ for the (abstract)
completion of $K$ at $\fP$ and $\iota_\fP$ for the natural embedding
$K\rightarrow K_\fP$. Let $K_p$ denote  the ring $\prod_{\fP\in S_p(K)}K_\fP$ endowed with
the product topology and the usual (continuous) $G$ action (see \eg~\cite[\S 2.3]{laser}). Thus the
diagonal embedding $\iota:=\prod_{\fP}\iota_\fP : K\rightarrow K_p$ is dense and
$G$-equivariant. We fix once and for all an algebraic closure
$\barbbQ_p$ of $\bbQ_p$ (equipped with the usual $p$-adic absolute
value $|\cdot|_p$), an embedding $j:\barbbQ\rightarrow\barbbQ_p$ and
a set $\tau_1,\ldots,\tau_d$ of left coset representatives for
$\Gal(\barbbQ/k)$ in $\Gal(\barbbQ/\bbQ)$. For each $i=1,\ldots,d$,
the embedding $j\circ\tau_i|_K:K\rightarrow\barbbQ_p$ extends to a
continuous embedding $K_{\fP_i}\rightarrow \barbbQ_p$ for a unique
prime ideal $\fP_i\in S_p(K)$, and we define
$\delta_i:K_p\rightarrow \barbbQ_p$ to be its composite with the
projection $K_p\rightarrow K_{\fP_i}$. (In general, the map
$i\mapsto \fP_i$ is not injective, nor surjective onto $S_p(K)$, but
the map $i\mapsto \fP_i\cap \cO_k$ is surjective onto $S_p(k)$.) For each
$\fP\in S_p(K)$, we write $U^1(K_\fP)$ for the group of principal
units of $K_\fP$ considered as a finitely generated $\bbZ_p$-module.
We write $U^1(K_p)$ for the group $\prod_{\fP\in S_p(K)}U^1(K_\fP)$
of `$p$-semilocal principal units of $K$' considered as a $\bbZ
G$-submobule of $K_p^\times$ and hence as a f.g.\ multiplicative
$\bbZ_p G$-module. ({\em Warning:} nevertheless, we shall often use
an {\em additive} notation for the $\bbZ_p G$-action on $U^1(K_p)$.)
It is clear that $|\delta_i(u)-1|_p<1$ for every $u\in U^1(K_p)$ and
each $i\in\{1,\ldots,d\}$ so that $\log_p(\delta_i(u))\in \barbbQ_p$
is given by the usual logarithmic series. The formula
$\lambda_{i,p}(u):=\sum_{g\in G}\log_p(\delta_i(gu))g\inv$ then
defines a $\bbZ_p G$-linear map
$\lambda_{i,p}:U^1(K_p)\rightarrow\barbbQ_p G$ and letting $i$ vary
we get a unique $\bbZ_p G$-linear `regulator' map $R_p$ from the
exterior power $\WUoKp$ to $\barbbQ_p G$ such that
$R_p(u_1\wedge\ldots\wedge u_d)=\det(\lambda_{i,p}(u_l))_{i,l=1}^d$.
(The dependence on $j$ of $\delta_i$, $\lambda_{i,p}$ and $R_p$ will
be denoted by a superscript `$j$' where necessary.) We can now
define a map
\beq\label{eq: def of sKkS}
\begin{array}{rcrcl}
 \sKkS &:& \WUoKp &\longrightarrow& \bbQ_p G^- \\
    & &    &                    \\
    & & \theta &\longmapsto    & j(a^{-,\ast}_{K/k,S})R^j_p(\theta)
\end{array}
\eeq
Some explanations
are in order. First, $x\mapsto x^\ast$ is the unique $\bbC$-linear
involution of $\bbC G$ sending $g$ to $g\inv$ for all $g\in G$.
Since $a^{-}_{K/k,S}$ lies in $\barbbQ G^-$, so does
$a^{-,\ast}_{K/k,S}$ and we apply $j$ coefficientwise to get an
element of $\barbbQ_p G^-$. Multiplying the result by
$R^j_p(\theta)$ in $\barbbQ_p G$ gives $\sKkS(\theta)$ which is {\em
a priori} another element of $\barbbQ_p G^-$. However one can show
that it actually lies in $\bbQ_p G^-$ and, moreover, is {\em
independent of the choice of $j$} (see~\cite[Prop.~3.4]{twizo}
and~\cite[Prop.~5]{laser}). Although the map $\sKkS$ is {\em not}
independent of the choice and ordering of the $\tau_i$'s, the
dependence is simple and explicit (see~\cite[Rem.~2.6]{laser} for
more details).

It is clear from its construction that $\sKkS$ is $\bbZ_p G$-linear. This implies
in particular that it vanishes on $\WUoKp^+$, so one loses nothing by regarding it as a map $\WUoKp^-\rightarrow\bbQ_p G^-$. This was the point of view
of~\cite{laser} but
for the present purposes it is slightly more convenient to take the domain to be the whole of $\WUoKp$. In this context, the statement (and proof)
of Prop.~6 of {\em ibid.} give
\begin{prop}\label{prop:ker and im of sKks}\ \\
\begin{enumerate}
\item\label{part: ker and im of sKks1} $\ker(\fs_{K/k,S})=\WUoKp^+ +\left(\WUoKp\right)_{\rm tor}$ and the second summand is finite.
\item\label{part: ker and im of sKks2} $\im(\fs_{K/k,S})$ spans $\bbQ_pG^-$ over $\bbQ_p$.\ePf
\end{enumerate}
\end{prop}
From now on we denote $\im(\fs_{K/k,S})$ by $\SKkS$.

The `Integrality Conjecture' of~\cite{laser} is precisely the statement that $\SKkS\subset\bbZ_p G^-$.
We shall not deal with this conjecture {\em per se} in the present paper
because, in the cases we shall study, it can be subsumed into the stronger Congruence Conjecture.
%% Finally, if one enlarges the set $S$ then
%% $a^{-}_{K/k,S}$ -- hence also $\sKkS$ and $\SKkS$ -- change in a simple way that follows from equation~(\ref{eq: introduction}) above and is described in equation~(23) of~\cite{laser}.

\subsection{Rubin-Stark Elements and the Pairing $H_{K/k,n}$}\label{sec: RS + pairing H}
Now let $\bar{G}:=\Gal(K^+/k)\cong G/\{1,c\}$. The $\bbC$-linear extension of the restriction map $G\rightarrow\bar{G}$ has kernel $\bbC G^-$ and defines an isomorphism
$\bbC G^+\cong\bbC\bar{G}$ identifying
$\Theta^+_{K/k,S}(s)=\sum_{\chi\ {\rm even }}L_{K/k, S}(s,\chi)e_{\chi\inv,G}$
with the function $\Theta_{K^+/k,S}(s)=\sum_{\chi\in \hat{\bar{G}}}L_{K^+/k, S}(s,\chi)e_{\chi\inv,\bar{G}}$.
Rubin-Stark elements
are conjectural elements of a certain exterior power of the $S$-units of $K^+$ which are
supposedly associated with the $d$th derivative of the Taylor series of $\Theta_{K^+/k,S}(s)$ at $s=0$. Indeed, using the functional equation for
primitive $L$-functions, one can show (see \eg~\cite[Ch.~I, \S 3]{Tate}) that
$\ord_{s=0}L_{K^+/k, S}(s,\chi)\geq d$ for all $\chi\in\hat{\bar{G}}$. (For $\chi=\chi_0$ one needs the fact that $|S|\geq d+1$, by Hypothesis~\ref{subhyp: S contains...}.)
Thus
\[
\Theta_{K^+/k,S}(s)=\Theta^{(d)}_{K^+/k,S}(0)s^d+o(s^d)\ \ \ \mbox{as $s\rightarrow 0$}
\]
where $\Theta^{(d)}_{K^+/k,S}(0)$ denotes the element
$\sum_{\chi\in \hat{\bar{G}}}\frac{1}{d!}\left(\frac{d}{ds}\right)^d|_{s=0}L_{K^+/k, S}(s,\chi)e_{\chi\inv,\bar{G}}$ of $\bbC \bar{G}$. Let
$e_{S,d,\bar{G}}$ be the (possibly empty) sum of the idempotents $e_{\chi\inv,\bar{G}}\in \bbC \bar{G}$
over those $\chi\in\hat{\bar{G}}$ for which $\ord_{s=0}L_{K^+/k, S}(s,\chi)$ is exactly $d$.
We refer to eq.~(13) of~\cite{laser} for an explicit formula for $e_{S,d,\bar{G}}$ demonstrating that it actually lies in $\bbQ \bar{G}$.
An element $m$ of any $\bbQ \bar{G}$-module $M$ will be said to {\em `satisfy the eigenspace condition (w.r.t.\ $(S,d,\bar{G})$)'} iff it lies in
$e_{S,d,\bar{G}}M$, \ie\ $m=e_{S,d,\bar{G}}m$. It is not hard to see that $\Theta^{(d)}_{K^+/k,S}(0)$
lies in $\bbR \bar{G}$ and satisfies the eigenspace condition. In fact, $\bbR \bar{G}\Theta^{(d)}_{K^+/k,S}(0)=e_{S,d,\bar{G}}\bbR \bar{G}$.

Let us write $U_S(K^+)$ for the group of all $S$-units of $K^+$, namely those elements of $K^{+,\times}$ which are local units at each place of $K^+$ above a place of $k$
which is not in $S$. We consider $U_S(K^+)$ as a multiplicative $\bbZ \bar{G}$-module and the tensor product
$\bbQ U_S(K^+):=\bbQ\otimes_\bbZ U_S(K^+)$ and its exterior power
$\WUSKp$ as natural $\bbQ \bar{G}$-modules. ({\em Warning:} we shall sometimes use an additive notation for these.)
For each $i=1,\ldots d$ we define a $\bbZ \bar{G}$-linear
map $\lambda_{i}:U_S(K^+)\rightarrow\bbR \bar{G}$ by
setting $\lambda_{i}(\eps):=\sum_{g\in G}\log|\tau_i(g\eps)|g\inv$. This `extends'
$\bbQ $-linearly to a map $\bbQ U_S(K^+)\rightarrow \bbR \bar{G}$, also denoted $\lambda_i$, which in turn
gives rise to a unique $\bbQ \bar{G}$-linear regulator map
$R_{K^+/k}$ from $\WUSKp$ to $\bbR \bar{G}$ such that $R_{K^+/k}(x_1\wedge\ldots\wedge x_d)=
\det(\lambda_i(x_l))_{i,l=1}^d$.

We now define a {\em Rubin-Stark element for $K^+/k$ and $S$} to be any element
$\eta$ of $\WUSKp$ satisfying the eigenspace condition w.r.t.\ $(S,d,\bar{G})$ and such that
\beq\label{eq: RS element}
\Theta^{(d)}_{K^+/k,S}(0)=R_{K^+/k}(\eta)
\eeq
One cannot currently demonstrate the existence of any $\eta\in\WUSKp$ satisfying~(\ref{eq: RS element})
unless either $K^+$ is absolutely abelian or all the characters $\chi\in\hat{\bar{G}}$ satisfying
$\ord_{s=0}L_{K^+/k, S}(s,\chi)=d$ are of order $1$ or $2$. On the other hand,
certain special cases of Stark's conjectures for the extension $K^+/k$
are essentially equivalent to the existence of such an $\eta$ (see~\cite[Rem.~2.3]{laser}) and one can, if necessary,
ensure that it simultaneously satisfies the eigenspace condition simply by replacing it by $e_{S,d,\bar{G}}\eta$.
This makes $\eta$ {\em unique} once $\tau_1,\ldots,\tau_d$, and hence $R_{K^+/k}$, have been fixed (\eg\ by~\cite[Lemma 2.7]{Rubin}).
Henceforth we shall therefore refer to such an element as {\em the} Rubin-Stark element for $K^+/k$ and $S$ and denote it
$\eta_{K^+/k,S}$. It may be thought of as a higher-order generalisation of a cyclotomic unit (or number).
%% Indeed, suppose $d=1$, $k=\bbQ$ and
%% $K^+=\bbQ(\mu_f)^+$, for some integer $f\geq 2$.  Provided also that $p|f$, the set $S_f:=\{\infty\}\cup\{q:q|f\}$
%% satisfies the conditions on $S$ and the Rubin-Stark element of
%% ${\textstyle \bigwedge^1_{\bbQ \bar{G}}\bbQ U_{S_f}(K^+)}=\bbQ U_S(K^+)$
%% is precisely $-\half\otimes (1-\xi_f)(1-\xi_f\inv)$.

From now on we shall assume that
\beq\label{hyp: K contains mupnpo}
\mbox{$K$ contains $\mu_\pnpo$}
\eeq
where $n$ is the integer of the Introduction, assumed w.l.o.g.\ to be $\geq 0$.
Thus, for each $\fP\in S_p(K)$,  $\iota_\fP$ induces an isomorphism $\mu_\pnpo(K)\rightarrow\mu_\pnpo(K_\fP)$ and
the local Hilbert symbol $(\alpha,\beta)_{K_\fP,\pnpo}\in\mu_\pnpo(K_\fP)$ is defined for any
$\alpha,\beta\in K_\fP^\times$. (We shall use
the definition of the Hilbert symbol given in~\cite{NeukirchANT} rather
than~\cite{Serre} which reverses the order of $\alpha$ and $\beta$,
thus effectively inverting $(\alpha,\beta)_{K_\fP,\pnpo}$.) Given any $\eps\in U_S(K^+)$ and $u=(u_\fP)_\fP\in U^1(K_p)$ we
define $[\eps,u]_{K,n}\in\ZpnpoZ$ by
\beq\label{eq: def of square brackets}
[\eps,u]_{K,n}=\sum_{\fP\in S_p(K)}\Ind_n\left(\iota_\fP\inv(\iota_\fP(\eps),u_\fP)_{K_\fP,\pnpo}\right)
\eeq
where $\Ind_n:\mu_\pnpo(K)\rightarrow\ZpnpoZ$ is the isomorphism defined by $\xi_\pnpo^{\Ind_n(\zeta)}=\zeta$ for all $\zeta\in\mu_\pnpo(K)$.
The pairing $[\cdot,\cdot]_{K,n}:U_S(K^+)\times U^1(K_p)\rightarrow\ZpnpoZ$ is bilinear and
one checks ({\em cf}~\cite[eq.\ (18)]{laser})
% [** check number in in final version***]
that
\beq\label{eq: G-equivariance of pairing to ZpnpoZ}
[g\eps,gu]_{K,n}=\kappa_n(g)[\eps,u]_{K,n}\ \ \ \mbox{for all $\eps\in U_S(K^+)$, $u\in U^1(K_p)$ and $g\in G$}
\eeq
where, here and henceforth, we write $\kappa_n$ for the {\em cyclotomic character modulo $\pnpo$}. We shall regard $\kappa_n$ as a homomorphism
$\Gal(\barbbQ/\bbQ)\rightarrow\ZpnpoZst$ whose restriction to $\Gal(\barbbQ/k)$ factors through
$G$ by~(\ref{hyp: K contains mupnpo}) and is denoted by the same symbol.
Thus, whether $g$ lies in $\Gal(\barbbQ/\bbQ)$ or in $G$, we have, by definition $g(\xi_\pnpo)=\xi_\pnpo^{\kappa_n(g)}$.
Next we consider the pairing $[\cdot,\cdot]_{K,n,G}$ defined as follows
\displaymapdef{[\cdot,\cdot]_{K,n,G}}{U_S(K^+)\times U^1(K_p)}{(\ZpnpoZ)G}{(\eps,u)}{\sum_{g\in G}[\eps,gu]_{K,n}g\inv}
If $h$ lies in $\bar{G}$ and  $\tilde{h}$ is any lift of $h$ in
$G$, then a short calculation using~(\ref{eq: G-equivariance of pairing to ZpnpoZ}) shows that
\beq\label{eq: G-equivariance of pairing to ZpnpoZG}
[h\eps,u]_{K,n,G}=\kappa_n(\tilde{h})\tilde{h}\inv[\eps,u]_{K,n,G}
\eeq
for any $\eps\in U_S(K^+)$ and $u\in U^1(K_p)$. Taking $h=1$, $\tilde{h}=c$ gives $[\eps,u]_{K,n,G}=-c[\eps,u]_{K,n,G}$ in $(\ZpnpoZ)G$. In other words,
$[\cdot,\cdot]_{K,n,G}$ takes values in $(\ZpnpoZ)G^-$. Let us denote by $\kappa_n^\ast$ the unique ring homomorphism
from $\bbZ \bar{G}$ to $(\ZpnpoZ)G^-$ which sends $h\in\bar{G}$ to
$\bar{2}\inv(\kappa_n(\tilde{h}_1)\tilde{h}_1\inv+\kappa_n(\tilde{h}_2)\tilde{h}_2\inv)$, where $\tilde{h}_1$ and $\tilde{h}_2=c\tilde{h}_1$ are the two lifts of $h$ to $G$. Then
equation~(\ref{eq: G-equivariance of pairing to ZpnpoZG}) shows that the pairing $[\cdot,\cdot]_{K,n,G}$ is $\kappa_n^\ast$-semilinear in the first variable.
% in the sense that, for any
% $x\in \bbZ \bar{G}$ we have
% $[x\eps,u]_{K,n,G}=\kappa_n^\ast(x)[\eps,u]_{K,n,G}$ for all $\eps$ and $u$.
On the other hand,
it follows from its definition that $[\cdot,\cdot]_{K,n,G}$ is  $\bbZ G$-linear, hence $\bbZ_p G$-linear, in the second variable. Consequently, we obtain a unique,
well-defined pairing $\cH_{K/k,n}:\bigwedge_{\bbZ\bar{G}}^d U_S(K^+)\times \WUoKp\rightarrow (\ZpnpoZ)G^-$ satisfying
\[
\cH_{K/k,n}(\eps_1\wedge\ldots\wedge\eps_d,u_1\wedge\ldots\wedge u_d)=\det([\eps_i,u_t]_{K,n,G})_{i,t=1}^d
\]
for any $\eps_1,\ldots,\eps_d\in U_S(K^+)$ and $u_1,\ldots,u_d\in U^1(K_p)$.
By construction, $\cH_{K/k,n}$ is
$\kappa_n^\ast$-semilinear in the first variable and $\bbZ_p G$-linear in the second and the latter implies
\beq\label{eq: H zero on WUoKp+}
\cH_{K/k,n}(\eta,\theta)=0\ \ \ \mbox{for all $\eta\in\bigwedge_{\bbZ\bar{G}}^d U_S(K^+)$ and $\theta\in\WUoKp^+$}
\eeq
(So, for fixed $\eta$ the map
$
\cH_{K/k,n}(\eta,\cdot):\WUoKp\rightarrow (\ZpnpoZ)G^-
$
factors through the projection on $\WUoKp^-$, just as $\sKkS$ does.) Finally, we can
`extend' $\cH_{K/k,n}$ in an obvious way so that the first variable lies in the tensor product $\Zbp\otimes_\bbZ\bigwedge_{\bbZ\bar{G}}^d U_S(K^+)$, where
$\Zbp$ denotes the subring $\{a/b\in\bbQ:p\ndiv b\}$ of $\bbQ$.

We now explain briefly a further `extension' of the pairing $\cH_{K/k,n}$ which is necessary to state the Congruence Conjecture properly but
-- for reasons that will become clear later -- has only a limited
importance for the computations of this paper. The reader may refer to~\cite{laser} for the details. Denote by
$\alpha_S$ the natural map $\bigwedge_{\bbZ\bar{G}}^d U_S(K^+)\rightarrow\bigwedge_{\bbQ\bar{G}}^d \bbQ U_S(K^+)$.
Following Rubin, we defined in {\em loc.\ cit.}, \S~2.2,
a $\bbZ \bar{G}$-lattice $\Lambda_{0,S}(K^+/k)$ in $\bigwedge_{\bbQ\bar{G}}^d \bbQ U_S(K^+)$ which contains the image of  $\alpha_S$ with finite index.
In {\em loc.\ cit.}, \S~2.3
we defined a pairing $H_{K/k,n}: \Zbp\Lambda_{0,S}(K^+/k)\times \WUoKp\rightarrow (\ZpnpoZ)G^-$ which has the following property.
If $1\otimes\alpha_S$ denotes the $\Zbp$-linearly extension of $\alpha_S$ to
$\Zbp\otimes\bigwedge^d_{\bbZ \bar{G}}U_{S}(K^+)$, then for any $\theta\in\WUoKp$ there is a commuting diagram
\beq\label{diag: H factors through H}
\xymatrix{
{\textstyle \Zbp\otimes\bigwedge^d_{\bbZ \bar{G}}U_{S}(K^+)}\ar[dd]_{1\otimes\alpha_S}\ar[drrrr]^{\cH_{K/k,n}(\cdot,\theta)}&&&&\\
&&&&(\bbZ/\pnpo\bbZ) G^-\\
\bbZ_{(p)}\Lambda_{0,S}(K^+/k)\ar[urrrr]_{H_{K/k,n}(\cdot,\theta)}&&&&
}
\eeq
This follows easily from~\cite[eq.~(20)]{laser}.\vspace{1ex}\\
\rem\label{rem: on the diagram}
In fact,
the vertical map above is an isomorphism whenever $p\ndiv|\bar{G}|$. (See~\cite[Remark~2.4]{laser}.) If $p||G|$, then both the kernel (namely the torsion in
$\Zbp\otimes\bigwedge^d_{\bbZ \bar{G}}U_{S}(K^+)$)
and the cokernel may be non-trivial, though finite.
As far as the present paper is concerned, the main consequence of~(\ref{diag: H factors through H}) is simply
that $\cH_{K/k,n}(\cdot,\theta)$ vanishes on the kernel of
$1\otimes\alpha_S$, for all $\theta$.

\subsection{Statement of the Conjecture}
With the above hypotheses the Congruence Conjecture ($CC$) of~\cite{laser} may be stated as follows.
\begin{conjalt}{$CC(K/k,S,p,n)$}\label{conj:CC}
The Rubin-Stark element $\eta_{K^+/k,S}$ exists and lies in $\bbZ_{(p)}\Lambda_{0,S}(K^+/k)$. Furthermore,
%% for each $\theta\in\WUoKp$ we have
if $\theta\in\WUoKp$ then $\sKkS(\theta)$ lies in $\bbZ_p G^-$ and satisfies the following congruence modulo $\pnpo$
\beq\label{eq: the congruence}
\overline{\sKkS(\theta)}=\kappa_n(\tau_1\ldots\tau_d)H_{K/k,n}(\eta_{K^+/k,S},\theta)\ \ \ \mbox{in $(\ZpnpoZ) G^-$.}
\eeq
\end{conjalt}
\rem\ The choice of $\tau_1,\ldots,\tau_d$ affects both $\sKkS$ and $\eta_{K^+/k,S}$ but not the validity of
$CC(K/k,S,p,n)$ thanks to the `normalising factor' $\kappa_n(\tau_1\ldots\tau_d)$ in~(\ref{eq: the congruence}).\vertsp\\
\rem\label{rem: changing K, S, n}\ The conjecture behaves well under changing $K$, $S$ and $n$.
More precisely, it is shown in~\cite[\S 5]{laser} that $CC(K/k,S,p,n)$ implies $CC(F/k,S',p,n')$ for any $S'$ containing $S$, any $n'$ such that $n\geq n'\geq 0$ and any intermediate field
$F$, $K\supset F\supset k$ provided that the norm map from $\WUoKp^-$ to ${\textstyle \bigwedge^d_{\bbZ_p \Gal(F/k)}U^1(F_p)^-}$ is surjective. This holds, for instance, if $K/F$
is at most tamely ramified at primes in $S_p(F)$.\vertsp\\
\rem\label{rem: $IC$} As already noted, the $CC$ includes the statement $\SKkS\subset\bbZ_p G^-$ \ie\ the Integrality Conjecture ($IC$). This was treated separately
in~\cite{laser} since it does not require $\mu_p\subset K$. Section~4 of {\em loc.\ cit.} contains a survey of evidence for both conjectures. The $IC$ is known in many cases where
the $CC$ is not, \eg\ when $p$ splits completely in $k$ (with a technical condition), when $p$ is unramified in $K$ or when $p\ndiv |G|$.\vertsp\\
\rem\label{rem: Rubin and the $CC$} It is shown in~\cite[Rem.~2.3]{laser} that Conjecture~$B'$ of~\cite{Rubin} implies the existence of $\eta_{K^+/k,S}$ and that it lies in
$\half\Lambda_{0,S}(K^+/k)$, hence it implies the first statement of the $CC$. However, even in situations where
$\eta_{K^+/k,S}$ is known as an explicit element of  $\half\Lambda_{0,S}(K^+/k)$ (for instance, if $k=K^+$)
the congruence~(\ref{eq: the congruence}) can still be elusive. If $\theta\in\WUoKp^+$ then~(\ref{eq: the congruence}) clearly holds trivially (\ie\ as $0=0$) and
the same thing happens in a couple of more interesting cases mentioned in~\cite[\S 4]{laser}. Apart from these, the full $CC$
is unknown whenever $K$ is not abelian over $\bbQ$.

\section{Methods of Computation}\label{sec: methods of comp}
In this section we describe in detail the method we used to
numerically check the $CC$ for the $48$ examples listed in
Section~\ref{sec: results}.

\subsection{The Set-Up}\label{subsec: set-up}

We take the field $k$ to be a real quadratic field (so $d=2$), and
always take $S=S^1$ (so we drop it from the notation when possible). In
view of Remark~\ref{rem: changing K, S, n}, $CC(K/k,S^1,p,n)$
implies $CC(K/k,S,p,n)$ for all other admissible $S$. The prime $p$
will be small for computational reasons: large primes would lead to
extensions $K/k$ of too large a degree. Thus we shall always take $p
= 3$, $5$ or $7$. For the same reason, we shall usually take $n =
0$, except for a few examples with $n = 1$ and $p = 3$, which were
added for completeness. Since $d = 2 < p$, the latter examples
necessarily have $p||G|$. The question as to whether or not $p$
divides $|G|$ is of importance in the computation of $\bbZ_p
G$-generators of $\WUotwo$, as we shall see in Subsection~\ref{sec:
gens of exterior square}. All computations take place in the number
field $F$ which is defined to be the {\em normal closure of $K$ over
$\bbQ$} (within $\barbbQ$). They were performed using the PARI/GP
system \cite{pari}.

\subsection{Computation of $a_{K/k}^-$ and $\sKk(\theta)$}\label{subsec: comp of sKk}

Using the implementation in PARI/GP of the algorithm of \cite{D-T}, see also
\cite[Section 10.3]{Co}, we can compute arbitrarily good
approximations of the values at $s= 1$ of the $S^1$-truncated Artin
$L$-functions of odd irreducible characters of $G$, and thus deduce
arbitrarily good approximations of $a_{K/k}^-$ as an element of
$i^2\mathbb{R}G=\mathbb{R}G$, thanks to~(\ref{eq: defn of a}). In order to
compute $\sKk$ we must however determine $a_{K/k}^-$ exactly as an
element of $FG^-$
%% where $F$ is some sufficiently large number field
%% containing $\mu_{f(K)}$ and $\sqrt{d_k}$ (see Section~\ref{subsec:
%%  TheMapsKkS}) in which choose to work.
and to this end, we use the
\begin{prop}\label{prop: stuff on atilde}
Let $f(K)$ denote the positive generator of the ideal $\ff(K)\cap\bbZ$. Set $\delta=1$ if $(p,f(K))=1$ and
$\delta=0$ otherwise, and let
\beq\label{eq: def of atildeminus}
\tilde{a}_{K/k}^-:=p^\delta |\mu(K)|\sqrt{d_k}N\ff(K)a_{K/k}^-=-p^\delta |\mu(K)|\sqrt{d_k}N\ff(K)\pi^{-2}\Theta^-_{K/k}(1)
\eeq
The coefficients of $\tilde{a}_{K/k}^-$ are algebraic
integers of $F\cap \bbQ(\mu_{f(K)})$ and are stable (as a set) under the
action of $\Gal(\barbbQ/\bbQ)$.
\end{prop}
The proof uses the following group-theoretic lemma whose (simple) verification is left to the reader.
\begin{lemma}
Let $\cal G$ be a group and $\cal H$ a subgroup of finite index in $\cal G$, and let ${\rm Ver}$ denote
the transfer homomorphism from ${\cal G}^{\rm ab}={\cal G}/{\cal G}'$ to ${\cal H}^{\rm ab}={\cal H}/{\cal H}'$. Suppose
$\cal J$ is a normal subgroup of $\cal H$ containing ${\cal H}'$ and write
$\tilde{\cal J}$ for the largest normal subgroup of $\cal G$ contained in $\cal J$, \ie $\tilde{\cal J}=\bigcap_g g{\cal J}g\inv$
where $g$ runs through ${\cal G}$ (or, indeed, through
a set of left-coset representatives for ${\cal H}$ in ${\cal G}$). Then $\tilde{\cal J}$ is contained in the kernel
of the composite homomorphism
\[
{\cal G}\longrightarrow{\cal G}^{\rm ab}\stackrel{\rm Ver}{\longrightarrow}{\cal H}^{\rm ab}\longrightarrow {\cal H}/{\cal J}
\]
\ePf
\end{lemma}
{\sc Proof of Proposition~\ref{prop: stuff on atilde}} Let $\Phi_{K/k}(s)$ be the function defined
in~\cite[eq.~(9)]{twizo}. It follows
from~\cite[eq.~(8)]{laser} (dropping $e^-$, since $k\neq \bbQ$) that
$\tilde{a}_{K/k}^-=(\prod(N\fp-\sigma_\fp\inv))|\mu(K)|d_k
N\ff(K)\Phi_{K/k}(0)$ where the product runs over the set of all
primes $\fp\in S_p(k)$ not dividing $\ff(K)$. (Since $K$ contains
$\mu_p$ and $[k:\bbQ]=2$, it is easy to see that either this set is
empty -- so $\delta=0$ -- or $p$ ramifies in $k$ and this set
consists of the unique prime $\fp\in S_p(k)$ -- so that
$N\fp=p=p^\delta$.) Equation~(27) of~\cite{twizo} shows that
the coefficients of $|\mu(K)|d_k
N\ff(K)\Phi_{K/k}(0)$ are algebraic
integers of $\bbQ(\mu_{f(K)})$, hence so are those of $\tilde{a}_{K/k}^-$.
%% The latter lie in $\bbR$ hence in $\bbQ(\mu_{f(K)})^+$.
It remains to show that they are $\Gal(\barbbQ/\bbQ)$-stable and lie in $F$.
Consider the automorphism of $\barbbQ G$ obtained
by applying some $\alpha\in\Gal(\barbbQ/\bbQ)$ to coefficients. It
was shown in~\cite[Prop.~3.2]{twizo} that this has the same effect on
$\Phi_{K/k}(0)$ as multiplying it by ${\cal V}_K(\alpha)$ where ${{\cal V}_K}$ is the composite homomorphism
\[
\Gal(\barbbQ/\bbQ)\longrightarrow\Gal(\barbbQ/\bbQ)^{\rm ab}\stackrel{\rm Ver}{\longrightarrow}
\Gal(\barbbQ/k)^{\rm ab}\longrightarrow G
\]
The same is therefore true of $\tilde{a}_{K/k}^-$, hence its coefficients are $\Gal(\barbbQ/\bbQ)$-stable. Moreover, they are fixed
by $\Gal(\barbbQ/F)$ because  the latter is contained in $\ker{{\cal V}_K}$, as follows from the Lemma. (Take ${\cal G}=\Gal(\barbbQ/\bbQ)$,
${\cal H}=\Gal(\barbbQ/k)$ and ${\cal J}=\Gal(\barbbQ/K)$, so that $\tilde {\cal J}=\Gal(\barbbQ/F)$).\ePf
\noindent\rem\ If we define $\tilde{a}_{K/k}^-$
to be the second member in~(\ref{eq: def of atildeminus}) with $p^\delta$ replaced by $\prod_{\fq\in S\setminus (S_{\rm ram}\cup S_\infty)}N\fq$, then
both the statement of the Proposition and its proof go through essentially unchanged for any $d>1$ and $S\supset S^1$. Since the coefficients
of $\tilde{a}_{K/k}^-$ also lie in $i^d\bbR$ in general, for $d=2$ they must actually lie in
$F\cap \bbQ(\mu_{f(K)})^+$.\vertsp\\
Let $\tilde{a}_{K/k, \sigma}^-$ denote the coefficient of $\sigma \in G$
in $\tilde{a}_{K/k}^-$. Having computed $a_{K/k}^-$ to high accuracy
in $\bbR G$ as described above, we obtain good real approximations
to the values $\tilde{a}_{K/k, \sigma}^-$ for $\sigma\in G$ and
hence to the coefficients of the polynomial $\prod_{\sigma \in G}(X
- \tilde{a}_{K/k, \sigma}^-)$. But Proposition~\ref{prop: stuff on
atilde} implies that this polynomial lies in $\bbZ[X]$, so we may
recover it exactly. By recognising the $\tilde{a}_{K/k, \sigma}^-$
among its roots in $F$ (embedded in $\bbC$), we then obtain
$\tilde{a}_{K/k}^-$ as an element of $\cO_F[G]$ and dividing by
$p^\delta |\mu(K)|\sqrt{d_k}N\ff(K)\in F^\times$ gives $a_{K/k}^-$
as an element of $FG$.

We now explain how to compute $\sKk(\theta)\in\bbQ_p G$ (for  $\theta\in\WUoKp$) to any predetermined ($p$-adic) accuracy. We shall need only
the case $\theta=u_1\wedge u_2$ with $u_1,u_2\in U^1(K_p)$ (which suffices anyway, by linearity).  For any integer $N\geq 1$
we write the power series $\log(1+X)$ as $\ell_N(X)+r_N(X)$  where
$\ell_N(X):=\sum_{t=1}^{N-1}(-1)^{t-1}X^t/t\in\bbQ[X]$ and
$r_N(X):=\sum_{t=N}^\infty(-1)^{t-1}X^t/t\in\bbQ[[X]]$.
For $i=1,2$ and any $u\in U^1(K_p)$
we define elements $\lambda_{i,p,N}(u)=\lambda_{i,p,N}^j(u)$ and $\rho_{i,p,N}(u)=\rho_{i,p,N}^j(u)$ of $\barbbQ_p G$ by
\[
\lambda_{i,p,N}(u):=\sum_{g\in G}\ell_N(\delta^j_i(g(u-1)))g\inv\ \ \ \mbox{and}\ \ \
\rho_{i,p,N}(u):=\sum_{g\in G}r_N(\delta^j_i(g(u-1)))g\inv
\]
so that $\lambda_{i,p}(u)=\lambda_{i,p,N}(u)+\rho_{i,p,N}(u)$. It follows easily that
$\lambda_{i,p}(u_l)=\lim_{N\rightarrow\infty}\lambda_{i,p,N}(u_l)$ for any $i,l\in\{1,2\}$
and consequently that
\[
\sKk(u_1\wedge u_2)=\lim_{N\rightarrow\infty}j(a_{K/k}^{-,\ast})\det(\lambda_{i,p,N}^j(u_l))_{i,l=1}^2
\]
The convergence in $\barbbQ_p G$ implied in each of these
limits is coefficientwise, w.r.t.\ the absolute value $|\cdot|_p$ on $\barbbQ_p$. The
next result gives us the explicit control we require on the rate of convergence  in the second limit.
First, we impose a $p$-adic norm $\|\cdot\|_p$ on the $\barbbQ_p$-algebra $\barbbQ_p G$ by setting
\[
\|a\|_p=\max\{|a_g|_p:g\in G\}\ \ \ \mbox{where $a=\sum_{g\in G}a_g g\in\barbbQ_p G$}
\]
It is easy to check that $\|x+y\|_p\leq\max\{\|x\|_p,\|y\|_p\}$ and $\|x.y\|_p\leq\|x\|_p.\|y\|_p$.
%% and $\|xa\|_p=|x|_p\|a\|_p$ for any $x\in\barbbQ_p$.
We define a rational number $m_{K/k}$ by $p^{m_{K/k}}=\|j(a_{K/k}^-)\|_p=\|j(a_{K/k}^{-,\ast})\|_p$. The coefficients
of  $a_{K/k}^-$ now being known as elements of a number field $F$,
we may calculate $m_{K/k}$ from their valuations at the prime ideal in $S_p(F)$ determined by $j$.
(Notice, however, that Prop.~\ref{prop: stuff on atilde} gives an {\em a priori} upper
bound for $m_{K/k}$ and also shows it to be independent of our choice of $j$.)
%% We set $\hat{m}_{K/k}=\max\{m_{K/k},0\}$ (which is typically equal to $m_{K/k}$).

Next, for $i=1,2$ we set $e_i=e_{\fP_i}(K/\bbQ)$ (recall that $\fP_i$ is  the element of $S_p(K)$ determined by $j\tau_i$) and
$h_i(x)=(\log(x)/\log(p))-(x/e_i)$ for any real number $x>0 $. Thus the function $h_i$ decreases monotonically to $-\infty$ on $[e_i/\log(p),\infty)$. For $i=1,2$,
we let $b_i$ be the smallest integer $b$ such that $p^b(p-1)\geq e_i$. In our examples, $b_i$ ranges from $0$ to $2$.
Finally, we write $\epsilon$ for the transposition $(1,2)\in \Sigma_2$.
\begin{prop}\label{prop: convergence of sKk}
Suppose that a  positive integer $M$ is given. Then for any integer $N>\max\{e_1,e_2\}/\log(p)$ satisfying the inequalities
\beq\label{eq: conditions for convergence of sKk}
h_{\epsilon(i)}(N)\leq -(M+m_{K/k}+(b_i-(p^{b_i}/e_i)))\ \ \ \mbox{for $i=1$ and $2$}
\eeq
we have
\beq\label{eq: estimate of convergence of sKk}
\|\sKk(u_1\wedge u_2)-j(a_{K/k}^{-,\ast})\det(\lambda_{i,p,N}^j(u_l))_{i,l=1}^2\|_p\leq p^{-M}\ \ \ \mbox{for all $u_1,u_2\in U^1(K_p)$}
\eeq
\end{prop}
\bPf\ If $t\in\bbZ_{\geq 1}$ and $i=1$ or $2$ then for any $g\in G$ and $u\in U^1(K_p)$ we clearly have
\[
|(\delta^j_i(g(u-1)))^t/t|_p\leq p^{-t/e_i}|t|_p\inv
\]
As $t$ varies, the R.H.S.\ of this inequality attains an absolute maximum of $p^{b_i-(p^{b_i}/e_i)}$ (at $t=p^{b_i}$) and, on the other hand, is always at most
$p^{h_i(t)}$. We deduce that for any $i$ and $u$, we have $\|\lambda_{i,p,N}(u)\|_p\leq p^{b_i-(p^{b_i}/e_i)}$
and $\|\rho_{i,p,N}(u)\|_p\leq p^{b_i-(p^{b_i}/e_i)}$ for {\em every}
positive integer $N$ and also
$\|\rho_{i,p,N}(u)\|_p\leq p^{h_i(N)}$ provided $N>e_i/\log(p)$. Therefore, writing
$\sKk(u_1\wedge u_2)$ as $j(a_{K/k}^{-,\ast})\det(\lambda_{i,p,N}^j(u_l)+\rho_{i,p,N}^j(u_l))_{i,l=1}^2$ and expanding the determinant, we find that for any
$N>\max\{e_1,e_2\}/\log(p)$ satisfying the inequalities~(\ref{eq: conditions for convergence of sKk}), we have
\[
\hspace*{4 em}\mbox{(L.H.S.\ of~(\ref{eq: estimate of convergence of sKk}))}
 \leq p^{m_{K/k}}\max\{p^{h_2(N)+b_1-(p^{b_1}/e_1)}, p^{h_1(N)+b_2-(p^{b_2}/e_2)}\}\leq p^{-M} \hspace{4 em}\Box
%% &\leq&\max\{p^{-M}, p^{\hat{m}_{K/k}+h_1(N)+h_2(N)}\}\\
%% &\leq &p^{-M}\max\{1, p^{-M+((p^{b_1}/e_1)-b_1)+((p^{b_2}/e_2)-b_2)}\}
\]
%% But the definition of $b_i$ shows that $(p^{b_1}/e_1)-b_1$ is at most $1/e_i$ (when $b_i=0$) which is in turn at most $1/(p-1)\leq 1/2$, since $\mu_p\subset K$ and $p\geq 3$. Therefore
%% $-M+((p^{b_1}/e_1)-b_1)+((p^{b_2}/e_2)-b_2)\leq 1-M\leq 0$ and the result follows.

\rem\label{rem: after convergence} \\
(i) The two inequalities~(\ref{eq: conditions for convergence of sKk}) coincide whenever $e_1=e_2$ and in particular, whenever $p$ does not split in $k$.
\vspace{1ex}\\
(ii) In our computations of $\sKk(u_1\wedge u_2)$, the elements
$u_1$, $u_2$ will always be `global' by which we shall mean that
$u_l=\iota(v_l)$ for $l=1,2$ where $v_l\in K^\times$ satisfies
$\ord_\fP(v_l-1)\geq 1$ for each $\fP\in S_p(K)$. (In fact, the $v_l$ will be constructed to lie in $\cO_K$.) Thus, for $i=1,2$
we can write $\lambda_{i,p,N}^j(u_l)$ as $j(\sum_{g\in G}(\tau_i
g(x_{l,N}))g\inv)$ where $x_{l,N}=\ell_N(v_l-1)$ lies in $K$ for
$l=1,2$, and so
\beq\label{eq: local to global, sort of!}
j(a_{K/k}^{-,\ast})\det(\lambda_{i,p,N}^j(u_l))_{i,l=1}^2=j\left(a_{K/k}^{-,\ast}\det({\textstyle
\sum_{g\in G}}\tau_i g(x_{l,N})g\inv)_{i,l=1}^2\right) \eeq
It follows from Proposition~\ref{prop: stuff on atilde} that the quantity inside the large parentheses on the R.H.S. of~(\ref{eq: local to
global, sort of!}) has coefficients in $F$. In fact, however, they lie in $\bbQ$. (Hints for a proof of this fact
are given in Rem.~3.3(i) and Props.~3.3 and~3.4 of~\cite{twizo} noting that $a_{K/k}^{-,\ast}=\sqrt{d_k}\Phi_{K/k}(0)^\ast$,
by~\cite[eq.~(8)]{laser}.) We may therefore drop the
`$j$' on the R.H.S. of~(\ref{eq: local to global, sort of!}) and substitute it into~(\ref{eq: estimate of convergence of sKk}).
%% The result is that when $u_1$ and $u_2$ are global,
%% Proposition~\ref{prop: convergence of sKk} provides, for each $M>1$,
%% an explicit construction of an element ${\mathfrak s}(u_1, u_2;M)$
%% of $\bbQ G^-$ such that $\sKk(u_1\wedge u_2)-{\mathfrak s}(u_1,
%% u_2;M)\in p^M\bbZ_pG^-$.

\subsection{Generators of $U^1(K_p)$}\label{subsec: generators of U1}
Both sides of the congruence~(\ref{eq: the congruence}) are $\bbZ_p G$-linear in $\theta$
so it suffices to test it on a set of $\theta$'s generating $\WUotwo$ over $\bbZ_p G$,
preferably few in number since the R.H.S. is particularly
computationally intensive (see Remark~\ref{rem: factor}).
First we explain our construction of a set $V$ of $\bbZ_p G$-generators for $U^1(K_p)$. This  is summarised in the two following
Propositions which hold
for any abelian extension of number fields $K/k$ and any
prime $p$.  For the rest of this subsection, we therefore drop the assumption $[k:\bbQ]=2$ and return to the notations of
Subsection~\ref{subsec: TheMapsKkS}. In addition, we shall write
$S_p(k)$ as $\{\fp(1),\ldots,\fp(t)\}$ where $t=|S_p(k)|\leq r$ and
$\fP(i,1),\ldots,\fP(i,h_i)$ for the distinct primes of $K$
dividing $\fp(i)$, for $i=1,\ldots,t$. For each pair $(i,j)$ with
$i=1,\ldots,t$ and $j=1,\ldots,h_i$ we shall abbreviate the
completion $K_{\fP(i,j)}$ to $\hat{K}_{i,j}$ and the embedding  $\iota_{\fP(i,j)}:K\rightarrow \hat{K}$  to $\iota_{i,j}$ so that $\iota=\prod_{i,j}\iota_{i,j}$ embeds $K$ in
$K_p=\prod_{i,j}\hat{K}_{i,j}$.
We write $\hat{\cO}_{i,j}$ for the ring of valuation integers of $\hat{K}_{i,j}$ and
$\hat{\fP}(i,j)$ for its maximal ideal. For any $l\geq 1$ we write $U_{i,j}^l$ for the
$l$th  term in the filtration of $\hat{\cO}_{i,j}^\times$, \ie\
$U_{i,j}^l=1+\hat{\fP}(i,j)^l$  considered as a finitely generated, multiplicative $\bbZ_p$-module. In particular $U^1(K_p)=\prod_{i,j}U^1_{i,j}\subset K_p^\times$.
We write $D_i$ for $D_{\fP(i,j)}(K/k)$ which depends only on $i$, since $K/k$ is abelian. The same is true for
$T_i:=T_{\fP(i,j)}(K/k)$,  $e'_i:=|T_i|=e_{\fP(i,j)}(K/k)$, $f_i:=|D_i/T_i|=f_{\fP(i,j)}(K/k)$ and for $\phi_i$ which we define to be
the Frobenius element at $\fP(i,j)$ considered as a generator of
the quotient group $D_i/T_i$. For each $i$ we also define a positive integer $l_i$ (independent of $j$) by
\[
l_i:=1+[pe_{\fP(i,j)}(K/\bbQ)/(p-1)]=1+[pe'_ie_{\fp(i)}(k/\bbQ)/(p-1)]
\]
It follows from the standard properties of $\log_p$ and $\exp_p$ as defined by the usual power series (see~\cite[Ch.~5]{Wash})
that $U_{i,j}^{l_i}$ is contained in $(U_{i,j}^{1})^p$ for all $i,j$.
Indeed, $u\in U_{i,j}^{l_i}$ implies $|u-1|_p<p^{-p/(p-1)}$ so that $|\frac{1}{p}\log_p(u)|_p=p|u-1|_p<p^{-1/(p-1)}$. Hence $v:=\exp_p(\frac{1}{p}\log_p(u))$ is a well-defined
element of $\hat{K}_{i,j}$
satisfying $v^p=u$ and $|v-1|_p<p^{-1/(p-1)}<1$, so that $v\in U_{i,j}^{1}$.
\begin{prop}\label{prop: gensofUone1}
In the above notation, suppose that for each $i$
we are given
\begin{enumerate}
\item \label{part: propgensofUone1} a subset $X_i=\{x_{i,1},\ldots, x_{i,n_i}\}$ of $\fP(i,1)$ such that the images of
$\iota_{i,1}(1+x_{i,1}),\ldots, \iota_{i,1}(1+x_{i,n_i})$ generate the quotient $U_{i,1}^1/U_{i,1}^{l_i}$ as a module over $\bbZ_p D_i$ and
\item \label{part: propgensofUone2} an element $a_i$ of $\cO_K$ such that for any $r=1,\ldots,t$ and $j=1,\ldots,h_r$ we have
\[
a_i\equiv\eitheror{1 \pmod{\fP(r,j)^{l_{r}}}}{if $r=i$ and $j=1$, and }{0 \pmod{\fP(r,j)^{l_{r}}}}{otherwise.}
\]
\end{enumerate}
Then the set $V:=\{\iota(1+a_ix_{i,s})\,:\, i=1,\ldots,t,\ s=1,\ldots,n_i\}$ generates $U^1(K_p)$ over $\bbZ_p G$.
\end{prop}
\bPf By Nakayama's Lemma, it suffices to show that $V$ generates $U^1(K_p)$ modulo $(U^1(K_p))^p$ and hence, by the preceding comments, that
the images of the elements $\iota(1+a_ix_{i,s})$ generate  $U^1(K_p)/\prod_{i,j}U^{l_i}_{i,j}$. As a $\bbZ_p G$-module this is the
product (over $i$) of the sub-modules $\prod_{j=1}^{h_i}\left(U^1_{i,j}/U^{l_i}_{i,j}\right)$ so that, by the definition of the $a_i$, it suffices to prove that for each $i$,
the latter is generated over $\bbZ_p G$
by the images of the elements $\iota_{i,1}(1+a_ix_{i,t})\times\ldots\times\iota_{i,h_i}(1+a_ix_{i,t})$ for $1\leq t\leq n_i$.
By the definition of the $a_i$ (again) and of $X_i$, these images lie in the subgroup $U^1_{i,1}/U^{l_i}_{i,1}$ and generate it over $\bbZ_p D_i$. Since
$\prod_{j}\left(U^1_{i,j}/U^{l_i}_{i,j}\right)$ is the direct product of $G$-translates of $U^1_{i,1}/U^{l_i}_{i,1}$, we are done.\ePf
\noindent To find a set $X_i$ as in part~\ref{part: propgensofUone1} of the statement of Proposition~\ref{prop: gensofUone1}, one could simply ensure that
the images of $\iota_{i,t}(1+x_{i,t})$ generate the finite module $U_{i,1}^1/U_{i,1}^{l_i}$ over $\bbZ_p$ rather than $\bbZ_p D_i$.
However, it is straightforward to construct a set that is generally smaller (for $f_i>1$), provided that the exact sequence
\beq\label{eq: splitting hypothesis}
\ses{1}{T_i}{D_i}{D_i/T_i}{1}
\eeq
splits. This is equivalent to the existence of a lift $\tilde{\phi}_i\in D_i$ of $\phi_i$ which is of order $f_i$. A sufficient condition is that $f_i$ be prime to
$e'_i$ or,
more generally, to the cardinality of $D_i^{f_i}$  (subgroup of $f_i$-th powers). Computational constraints on $[K:k]$ mean that
$D_i$ is a fairly small group in the examples considered, so it is perhaps not so surprising that the sequence~(\ref{eq: splitting hypothesis})
was found to split in {\em all} of them (for all $i$) without any pre-selection. Assuming that this occurs, we write $A_i$ for the subgroup (of order $f_i)$ of
$D_i$ generated by some $\tilde{\phi}_i$,
$N_i$ for $K^{A_i}$ and $\wp_i$ for the prime of $N_i$ below $\fP(i,1)$:
\[
\xymatrix{
K\ar@/_1pc/@{-}[d]_{A_i}\ar@{-}[d]&&\fP(i,1)\ar@{-}[d]\\
N_i\ar@{-}[d]&&\wp_i\ar@{-}[dd]\\
K^{D_i}\ar@{-}[d]\\
k&&\fp(i)
%\bbQ(f)\ar@/^2.5pc/@{-}[ddd]^{G(f)}\ar@{-}[ddd]\ar@/_1pc/@{-}[dl]_{\langle c\rangle}\ar@{-}[dl]\\
%\bbQ(f)^+\ar@{-}[ddr]\\
%\\
%&\bbQ
}
\]
Thus $N_i/K^{D_i}$ is totally ramified at $\wp_i$  (so $f_{\wp_i}(N_i/k)=1$). On the other hand,
$K/N_i$  is unramified at $\fP(i,1)$ and $A_i=\Gal(K/N_i)$ maps
isomorphically onto the Galois group of the residue field $\cO_K/\fP(i,1)$ over $\cO_{N_i}/\wp_i=\cO_k/\fp(i)$ (with $\tilde{\phi}_i$ acting by $N\fp(i)$-th powers).
It follows from the Normal Basis Theorem that $\cO_K/\fP(i,1)$ is freely generated over $(\cO_k/\fp(i))A_i$.
Moreover, a well-known criterion states that a
free generator is given by the class $\bar{\alpha}_i$ modulo $\fP(i,1)$ of $\alpha_i\in\cO_K$ if and only if
$\det(gh(\bar{\alpha_i}))_{g,h\in A_i}\neq 0$ in $\cO_K/\fP(i,1)$; in other words, iff
$\det(\alpha_i^{N\fp(i)^{a+b}})_{a,b=0}^{f_i-1}\nin\fP(i,1)$. Such
an $\alpha_i$ is easily found
by trial and error.
\begin{prop}\label{prop: gensofUone2} Suppose the sequence~(\ref{eq: splitting hypothesis}) splits for some $i\in\{1,\ldots, t\}$. Using the above notations,
choose any $\pi_i\in\wp_i\setminus\wp_i^2$ and a subset $Y_i\subset\cO_k$ whose images in $\cO_k/\fp(i)$ form a basis over
$\bbZ/p\bbZ$. Then the set $X_i:=\{\pi_i^a y \alpha_i\,:\, a=1,\ldots,l_i-1,\ y\in Y_i\}$ satisfies the requirements of part~\ref{part: propgensofUone1} of the statement of
Proposition~\ref{prop: gensofUone1} (and even with `over $\bbZ_p D_i$', replaced by `over $\bbZ_p A_i$').
\end{prop}
\bPf\ The definition of $\alpha_i$ ensures that the classes modulo $\hat{\fP}(i,1)$ of the $\iota_{i,1}(y \alpha_i)$ for $y\in Y_i$ freely generate
$\hat{\cO}_{i,1}/\hat{\fP}(i,1)$ over $(\bbZ/p\bbZ) A_i$. Now
$\iota_{i,1}(\pi_i)$ is a local uniformiser for $K_{i,1}$ so for each $a=1,\ldots,l_i-1$ there is a familiar isomorphism of (finite) $\bbZ_p$-modules
$\hat{\cO}_{i,1}/\hat{\fP}(i,1)\rightarrow U^{a}_{i,1}/U^{a+1}_{i,1}$ which sends the class of $x$ to the class of $1+\iota_{i,1}(\pi_i)^ax$. {\em Since  $\iota_{i,1}(\pi_i)$ is fixed by $A_i$}, this is a
$\bbZ_p A_i$-isomorphism and it follows that for each $a=1,\ldots,l_i-1$
the classes modulo $U^{a+1}_{i,1}$ of the $\iota_{i,1}(1+\pi_i^ay \alpha_i)$ for $y\in Y_i$ generate $U^{a}_{i,1}/U^{a+1}_{i,1}$ over $\bbZ_p A_i$.
The result follows easily from this.\ePf
\noindent For each example tested, we used Propositions~\ref{prop:  gensofUone1} and~\ref{prop: gensofUone2} to construct a set of $\bbZ_p G$-generators
for $U^1(K_p)$ which is denoted $V$ and has cardinality
$N:=\sum_{i=1}^t|X_i|=\sum_{i=1}^t(l_i-1)|Y_i|=\sum_{i=1}^t(l_i-1)f_{\fp(i)}(k/\bbQ)\leq\frac{pd}{p-1}\max\{e'_i:i=1,\ldots,t\}$.\vertsp\\
\rem\label{rem: modifying v's} By construction, each  $u\in V$  is of form $\iota(v)$ for some
$v\in\cO_K$ which is congruent to $1$ modulo each $\fP\in S_p(K)$. At certain points in the computations it can be helpful to
`perturb' one or several such $v$ as follows
\[
\mbox{$v\leadsto v':=v+p^{l+1}x$ for some $x\in\cO_K$ and $l\in\bbZ$, $l\geq 1$}
\]
Clearly, $\iota(v')\in U^1(K_p)$ and previous arguments involving $\log_p$ and $\exp_p$ can be adapted to show that
\[
\iota(v')\equiv\iota(v)\pmod{U^1(K_p)^{p^l}}
\]
(Indeed,  $\iota_{i,j}(v'/v)$ is the $p^l$th power of $\exp_p(p^{-l}\log_p(1+p^{l+1}\iota_{i,j}(x/v)))\in U^1(\hat{K}_{i,j})$, convergence being
assured by  the fact that $l+1>l+(1/(p-1))$ since $p>2$.) In particular,
Nakayama's Lemma implies that such perturbations do not effect the generation of $U^1(K_p)$ by $V$. For example, taking $l=1$,
one can modify the $v$'s corresponding to each $u\in V$ to
ensure that their coefficients with respect to a given $\bbZ$-basis of $\cO_K$ have absolute value at most $p^2/2$.

\subsection{Generators of $\WUotwo$}\label{sec: gens of exterior square}

Proposition~\ref{prop:ker and im of sKks}\ref{part: ker and im of sKks1} and Equation~(\ref{eq: H zero on WUoKp+}) show
that both sides of~(\ref{eq: the congruence}) vanish for $\theta\in\WUotwo^+$, so we only need a set of generators modulo this submodule. For the L.H.S.\ of~(\ref{eq: the congruence}),
Proposition~\ref{prop:ker and im of sKks}\ref{part: ker and im of sKks1} shows
that the same is true for $\theta\in\left(\WUotwo\right)_{\rm tor}$.
However for the R.H.S.\ we have only managed to prove this under the assumption $p\ndiv |G|$ (see
\cite[Proposition 8]{laser}), so we proceed as follows.
%% Furthermore, under the same assumption, Equation (33) of {\em ibid.} proves the existence of such a set consisting of just {\em one} element
%% which is of form $u_1\wedge\ldots\wedge u_d$. However,
%% the proof is non-constructive and
%% In all cases we first compute a set $V=\{v_1,\ldots,v_N\}$ of $\bbZ_p G$-generators for $U^1(K_p)$, aiming to minimise $N$. The
%% $\left(
%% \begin{array}{c}
%% N \\
%% 2
%% \end{array}\right)$
For the minority of examples considered where $p$ divides $|G|$,
we simply test~(\ref{eq: the congruence})
for all $\theta$ in the $\half N(N-1)$-element set
$W:=\{v_s\wedge v_r\,:\,1\leq s<r\leq N\}$, with $V = \{v_1, \dots, v_N\}$, which clearly generates all of $\WUotwo$ over $\bbZ_p G$.
For the rest of this subsection we will assume $p\ndiv |G|$ and describe a second procedure to construct a
subset $W'\subset W$  generating $\WUotwo$ modulo $\WUotwo^++\left(\WUotwo\right)_{\rm tor}$ and such that $|W'|$ is much smaller than $|W|$ (see below). By the above remarks, it
will then suffice to test~(\ref{eq: the congruence}) for all $\theta$ in $W'$.
A generic element of $W$ will be denoted $\theta$. Even for integers $M$
somewhat greater than $n+1$, the computation of
$\sKk(\theta)$ modulo $p^M\bbZ_p G^-$ is relatively quick
compared with that of $H_{K/k,n}(\eta_{K^+/k,S},\theta)$ in $(\ZpnpoZ) G^-$.
We turn this fact to our advantage by using $\sKk$ itself to determine  $W'$.
Indeed, it is obvious from Proposition~\ref{prop:ker and im of sKks}\ref{part: ker and im of sKks1}
that $W'$ will have the required property if and only if $\SKk$ equals the $\bbZ_p G^-$-submodule
$\langle\sKk(\theta): \theta\in W'\rangle_{\bbZ_p G^-}$ of $\bbQ_pG^-$.
%% is the
%% \beq\label{eq: required equality}
%% \langle\sKk(\theta): \theta\in W'\rangle_{\bbZ_p G^-}=
%% \langle\sKk(\theta): \theta\in W\rangle_{\bbZ_p G^-}(=\SKk)
%% \eeq
%% $\bbZ_p G^-$-submodule of $\bbQ_p G^-$ generated by $\{\sKk(\theta): \theta\in W\}$.
We construct such a $W'$ by means of an explicit isomorphism from $\bbQ_p G^-$ to a product of fields which we now describe.

Since $G$ is small, it is easy to compute
a set $R^-$ of representatives of the orbits of the odd, irreducible characters $\chi:G\rightarrow\barbbQ^\times$ under the action of $\Gal(\barbbQ/\bbQ)$.
The $\bbQ$-linear
extension of each such character $\chi$ defines a homomorphism from $\bbQ G^-$ to $F_\chi:=\bbQ(\chi)$ such that the product over $\chi\in R^-$ is a
ring isomorphism $X^-:\bbQ G^-\rightarrow\prod_{\chi\in R^-}F_\chi$.
Tensoring
over $\bbQ$ with $\bbQ_p$ we get the first isomorphism, $X_p^-$, below.
\beq\label{eq: the iso alpha}
\bbQ_p G^-\stackrel{X_p^-}{\longrightarrow}\prod_{\chi\in R^-}(F\chi\otimes\bbQ_p)\stackrel{Z_p^-}\longrightarrow \prod_{\chi\in R^-}\prod_{\fP\in S_p(F_\chi)}F_{\chi,\fP}
\eeq
The second, $Z_p^-$, is the product over $\chi\in R^-$ of the isomorphisms from $F_\chi\otimes\bbQ_p$ to the product of the completions of $F_\chi$
at primes above $p$, the latter taking
$a\otimes x$ to the vector $(x\iota_\fP(a))_\fP$.
Let us write the composite isomorphism $Z_p^-\circ X_p^-$ as $\alpha=\prod_{\chi}\prod_\fP\alpha_{\chi,\fP}$. We identify $\bbZ_p G^-$ with $((1-c)\bbZ G)\otimes\bbZ_p$ considered as a subring of $\bbQ G^-\otimes \bbQ_p$ which we are identifying
with $\bbQ_pG^-$. It is clear that $X_p^-$ sends $\bbZ_p G^-$
into $\prod_{\chi\in R^-}(\cO_\chi\otimes\bbZ_p)$ where $\cO_\chi:=\bbZ[\chi]=\cO_{F_\chi}$ and that the image surjects onto each component (since $p\neq 2$). For a
given $\chi\in R^-$, let  us write $e(\chi)$ for
the sum of the idempotents in $\barbbQ G$ belonging to the irreducible characters in the $\Gal(\barbbQ/\bbQ)$-orbit of $\chi$. It is easy to see that
$e(\chi)$ lies in $(1-c)|G|\inv\bbZ G$ inside $\bbQ G^-$ and
hence that $e(\chi)\otimes 1$ lies in
$\bbZ_p G^-$ (since $p\ndiv |G|$). The orthogonality relations imply  that $X_p^-(e(\chi)\otimes 1)$ has component $1$ at $\chi$ and $0$ elsewhere.
It follows that $\bbZ_p G^-$ is sent  isomorphically onto $\prod_{\chi\in R^-}(\cO_\chi\otimes\bbZ_p)$ by $X_p^-$. Hence $\alpha$ maps it
isomorphically onto the image of the latter under $Z_p^-$ which, by standard facts, is $\prod_{\chi\in R^-}\prod_{\fP\in S_p(F_\chi)}\cO_{\chi,\fP}$,
where $\cO_{\chi,\fP}$ denotes the ring of integers of $F_{\chi,\fP}$. The values of $\chi$ are roots of unity of order prime to $p$. It follows that $F_\chi/\bbQ$ is unramified at
$p$ so that each $\cO_{\chi,\fP}$ is a complete d.v.r. with maximal ideal $p\cO_{\chi,\fP}$.

Both $\SKk$ and $\bbZ_p G^-\sKk(\theta)$ (for any $\theta\in W$) are $\bbZ_pG^-$ submodules
of $\bbQ_pG^-$. Hence, for each pair $(\chi,\fP)$ and each $\theta\in W$ there exist $m(\chi,\fP)$ and $m(\chi,\fP;\theta)$ in $\bbZ\cup\{\infty\}$
such that (taking $p^\infty=0$):
\[
\mbox{
${\displaystyle
\alpha(\SKk)=\prod_{\chi\in R^-}\prod_{\fP\in S_p(F_\chi)}p^{m(\chi,\fP)}\cO_{\chi,\fP}
}$
\ \ and\ \
${\displaystyle
\alpha(\bbZ_p G^-\sKk(\theta))=\prod_{\chi\in R^-}\prod_{\fP\in S_p(F_\chi)}p^{m(\chi,\fP;\theta)}\cO_{\chi,\fP}
}$
}
\]
Of course, $m(\chi,\fP;\theta)$ is just $\ord_p(\alpha_{\chi,\fP}(\sKk(\theta)))$, while $m(\chi,\fP)<\infty$ by Proposition~\ref{prop:ker and im of sKks}\ref{part: ker and im of sKks2}
and $m(\chi,\fP)\geq 0$ since the Integrality Conjecture is known for $p\ndiv |G|$ by~\cite[Cor.~1]{laser}.
The properties of d.v.r.'s give the equivalence
\beq\label{eq: equivalence}
\SKk=\langle\sKk(\theta): \theta\in W\rangle_{\bbZ_p G^-}\Longleftrightarrow m(\chi,\fP)=\min\{m(\chi,\fP;\theta):\theta\in W\}\
\mbox{for all $(\chi,\fP)$}
\eeq
Since the first equality holds by construction of $W$, so must the second and in particular $m(\chi,\fP;\theta)\geq 0\ \forall \chi,\fP,\theta$.
For each pair $(\chi,\fP)$ we define $W_{\min}(\chi,\fP)$ to be the non-empty subset of $W$ on which  $m(\chi,\fP;\theta)$ attains its minimum.
By the above,
$W_{\min}(\chi,\fP)=\{\theta\in W:m(\chi,\fP;\theta)=m(\chi,\fP)\}$.

The construction of $W'$ begins by using Proposition~\ref{prop: convergence of sKk}
to compute an approximation to $\sKk(\theta)$ for each $\theta\in W$,
with a guaranteed $p$-adic precision of
$p^{-M}$ for a moderate value of $M\geq n+1$,
\eg\ $M=n+3$. Since each $\theta$ is already expressed as $u_1\wedge u_2$ for `global' elements $u_1$ and $u_2$ in the sense of Remark~\ref{rem: after convergence}(ii)),
the latter and Proposition~\ref{prop: convergence of sKk} naturally give rise to an
approximation in $\bbQ G^-$, which we somewhat abusively write as $\sKk(\theta;M)$, such that $\sKk(\theta)-\fs(\theta;M)\in p^M\bbZ_pG^-$.
Now, fixing  $\chi\in R^-$ and $\fP\in S_p(F_\chi)$,
we may compute the values
$\ord_p(\alpha_{\chi,\fP}(\fs(\theta;M)))$ one by one for each $\theta\in W$, since these are just $\ord_\fP(\chi(\fs(\theta;M)))$, by construction of $\alpha$. Suppose
$W_{<M}(\chi,\fP,M)$ is the subset of those $\theta$ in $W$ for which we find $\ord_p(\alpha_{\chi,\fP}(\fs(\theta;M)))<M$.
By the ultrametric inequality, if $\theta$ lies in $W_{<M}(\chi,\fP,M)$ then we must have $m(\chi,\fP;\theta)=\ord_p(\alpha_{\chi,\fP}(\fs(\theta;M)))$. Otherwise we know only
that $m(\chi,\fP;\theta)\geq M$. This means that if $W_{<M}(\chi,\fP,M)$ is non-empty, we may compute
%% $m(\chi,\fP)$ as the minimum value of $\ord_p(\alpha_{\chi,\fP}(\fs(\theta;M)))$ for
%% $\theta\inW_{<M}(\chi,\fP,M)$ and
$W_{\min}(\chi,\fP)$ as the subset of $W_{<M}(\chi,\fP,M)$ on which $\ord_p(\alpha_{\chi,\fP}(\fs(\theta;M)))$ attains its minimum,
and then pass to the next pair $(\chi,\fP)$. However, in a very small number
of examples we encountered a pair $(\chi,\fP)$ for which $W_{<M}(\chi,\fP,M)=\varnothing$ for the initial value of $M$. In this case we simply recalculated the $\fs(\theta;M)$'s
with a larger value of $M$ until $W_{<M}(\chi,\fP,M)\neq\varnothing$ for that pair and then continued with the increased value of $M$. This simple hit-and-miss procedure
terminated rapidly enough: in all cases we were able to determine $W_{\min}(\chi,\fP)$ for all pairs $(\chi,\fP)$ without ever taking $M>n+5$.

The equivalence obtained by replacing $W$ by $W'$ in~(\ref{eq: equivalence}) shows that  a subset $W'\subset W$
will have the required property iff $W'\cap W_{\min}(\chi,\fP)\neq\varnothing$ for all pairs
$(\chi,\fP)$. Picking an
element at random from each $W_{\min}(\chi,\fP)$ would give a subset $W'$ whose cardinality could not exceed
the number of pairs $(\chi,\fP)$ which in turn is at most ${\rm dim}_{\bbQ_p}\bbQ_pG^-=\half|G|$, by~(\ref{eq: the iso alpha}). This is already much smaller than $|W|$ in most cases.
In practice, however, there was a tendency for $\bigcap_{\chi,\fP}W_{\min}(\chi,\fP)$ to be non-empty, so we could simply take $W'=\{\theta_0\}$
for any $\theta_0$ in this intersection. This tendency is explained by the fact that, as a submodule
of finite index in $\bbZ_p G^-$ (which is a product of d.v.r.'s),  $\SKk$   is automatically free over $\bbZ_pG^-$ with one generator.
While there is no guarantee that
$W$ contains such a generator, it is not surprising that it often does.
In fact, if this failed for our initial choice of $W$, the best practical solution was simply to randomly  modify $W$
once or twice  until it did
(\eg\ by changing the elements $a_i$ and $x_{i,j}$ used in
Prop.~\ref{prop: gensofUone1} to construct $V$). We thus achieved
$|W'|=1$ in all cases without too much difficulty. \vertsp\\
\rem\ The procedure described above for $p\ndiv |G|$ determines the values $m(\chi,\fP)$ for all $(\chi,\fP)$ as a by-product.
However, they are also given
explicitly by `index formula'~(31) of~\cite{laser} (see also~(32) of {\em ibid.}). One should take `$\phi$'
to be the composition of $\chi$ with
any embedding $\barbbQ\rightarrow\barbbQ_p$ inducing $\fP$. Using this formula, one could in principle select the
initial value of $M$ to be greater than the maximum of the
$m(\chi,\fP)$'s, thereby ensuring that $W_{<M}(\chi,\fP,M)\neq\varnothing$ for all $(\chi,\fP)$.

\subsection{Computation of $\eta_{K^+/k,S^1}$}\label{subsec: compute rs-element}
We need to determine the Rubin-Stark element $\eta_{K^+/k,S^1}$, that is,
the unique (conjectural) element of ${\textstyle \bigwedge_{\bbQ \bar{G}}^2}\bbQ U_{S^1}(K^+)$ that
satifies the eigenspace condition w.r.t.\ $(S^1,d,\bar{G})$ and
Equation~(\ref{eq: RS element}). The first statement in the Congruence Conjecture tells us to expect
$\eta_{K^+/k,S}$ to lie in $\bbZ_{(p)}\Lambda_{0,S}(K^+/k)$.
(As noted in Remark~\ref{rem: Rubin and the $CC$}, this is also predicted by Conjecture~$B'$ of~\cite{Rubin}.)
Assuming this {\em and also $p\ndiv|\bar{G}|$}, we shall now sketch a proof that $\eta_{K^+/k, S^1}$ must in fact be of the rather more precise form
\beq\label{eq: what to expect}
\eta_{K^+/k, S^1}=(1\otimes \alpha_{S^1})\left(\frac{1}{a}\otimes(\eps_1\wedge\eps_2)\right)
%%\frac{1}{a}\otimes(\eps_1\wedge\eps_2)
\ \ \ \mbox{for some $\eps_1,\eps_2\in U_{S(p)}(K^+)$ and $a\in\bbZ$, $p\ndiv a$}
\eeq
where $S(p) := S_\infty \cup S_p$. First, by Remark~\ref{rem: on the diagram}, these hypotheses imply
$\eta_{K^+/k, S^1}=(1\otimes \alpha_{S^1})(\tilde{\eta})$ for a (unique)
$\tilde{\eta}\in \bbZ_{(p)}\otimes \bigwedge^d_{\bbZ \bar{G}}U_{S^1}(K^+)$. Identifying the latter module
with $\bigwedge^d_{\bbZ_{(p)} \bar{G}}\bbZ_{(p)}\otimes U_{S^1}(K^+)$, we may write $\tilde{\eta}=\sum_{i=1}^Nx_{1,i}\wedge x_{2,i}$ where
$x_{1,i},x_{2,i}\in \bbZ_{(p)}\otimes U_{S^1}(K^+)$. Writing $e_{S^1}$ for the idempotent
$e_{S^1,2,\bar{G}}\in|\bar{G}|\inv\bbZ \bar{G}\subset \bbZ_{(p)} \bar{G}$, the eigenspace condition gives
\[
\eta_{K^+/k, S^1}=e_{S^1}(1\otimes \alpha_{S^1})(\tilde{\eta})=(1\otimes \alpha_{S^1})(e_{S^1}\tilde{\eta})
=(1\otimes \alpha_{S^1})(\sum_{i=1}^Ne_{S^1}x_{1,i}
\wedge e_{S^1}x_{2,i})
\]
Now consider
$A_{S^1}:=e_{S^1}(\bbZ_{(p)}\otimes U_{S^1}(K^+))$ as a module over
the ring $e_{S^1}\bbZ_{(p)} \bar{G}$ which is a product of p.i.d.'s (again because $p\ndiv|\bar{G}|$). Since $p\neq 2$, $A_{S^1}$ is
$\bbZ$-torsionfree. Moreover $\bbQ A_{S^1}$ is free of rank $2$ over $e_{S^1}\bbQ \bar{G}$. (This follows from the definition of $e_{S^1,2,\bar{G}}$
and the fact that $\dim_\bbC(e_{\chi\inv,\bar{G}}\bbC U_{S^1}(K^+))=\ord_{s=0}L_{K^+/k,S^1}(s,\chi)$ for all $\chi\in\hat{\bar{G}}$.) It follows easily that
$A_{S^1}$ is free of rank $2$ over $e_{S^1}\bbZ_{(p)} \bar{G}$ and it is not hard to see that any pair of basis elements can be written
$\frac{1}{a_1}\otimes \alpha_1$, $\frac{1}{a_2}\otimes\alpha_2$ where $p\ndiv a_1, a_2\in\bbZ$ and
and $\alpha_1,\alpha_2$ lie in $(|\bar{G}|e_{S^1})U_{S^1}(K^+)$. Writing each $e_{S^1}x_{1,i}$ and
$e_{S^1}x_{2,i}$ in such a basis, we conclude that $\tilde{\eta}$ is a $\bbZ_{(p)}\bar{G}$-multiple of
$1\otimes(\alpha_1\wedge\alpha_2)$. Equation~(\ref{eq: what to expect}) will
clearly follow if we can show $(|\bar{G}|e_{S^1})U_{S^1}(K^+)\subset U_{S(p)}(K^+)$.
If $S^1=S(p)$, this is immediate. Otherwise $|S^1|> d+1$ and~\cite[eq. (13)]{laser}
shows that $N_{D_\fq}e_{S^1}=0$ where $N_{D_\fq}\in\bbZ\bar{G}$ is the norm element of
the decomposition subgroup $D_\fq\subset\bar{G}$ for any prime $\fq\in S^1\setminus S_\infty$. Thus every element of
$(|\bar{G}|e_{S^1})U_{S^1}(K^+)$ is killed by every such $N_{D_\fq}$ which implies that, in fact,
$(|\bar{G}|e_{S^1})U_{S^1}(K^+)\subset U_{S_\infty}(K^+)$ \ie\ in this case, we can actually take $\eps_1,\eps_2\in\cO_{K^+}^\times$
in~(\ref{eq: what to expect}).

Let us write $\cU(p)$ for $\bbQ\otimes{\textstyle \bigwedge_{\bbZ \bar{G}}^2} U_{S(p)}(K^+)$ considered as a $\bbQ\bar{G}$-submodule of
${\textstyle \bigwedge_{\bbQ \bar{G}}^2}\bbQ U_{S^1}(K^+)$. If $p\ndiv |\bar{G}|$, we have just shown that the Congruence
Conjecture implies~(\ref{eq: what to expect}) or, equivalently,
\beq\label{eq: what to expect 2}
\eta_{K^+/k, S^1}=
\frac{1}{a}\otimes(\eps_1\wedge\eps_2)\in e_{S^1}\cU(p)
\ \ \ \mbox{for some $\eps_1,\eps_2\in U_{S(p)}(K^+)$ and $a\in\bbZ$, $p\ndiv a$}
\eeq
If $p||\bar{G}|$ and we assume only that the Rubin-Stark element $\eta_{K^+/k,S^1}$ exists,
then similar but simpler arguments (replacing  $\Zbp$ by $\bbQ$)
still show that
\beq\label{eq: what to expect 3}
\eta_{K^+/k, S^1}=
\frac{1}{a}\otimes(\eps_1\wedge\eps_2)\in e_{S^1}\cU(p)
\ \ \ \mbox{for some $\eps_1,\eps_2\in U_{S(p)}(K^+)$}
\eeq
and also that
$e_{S^1}\cU(p)$ is free of rank $1$ over $e_{S^1}\bbQ\bar{G}$.
These observations motivate the following procedure for determining $\eta_{K^+/k, S^1}$
which is much simpler than the one used in \cite{Solomon-Roblot1} but still sufficient for present purposes.
First we compute an $e_{S^1}\mathbb{Q}\bar{G}$-generator of $e_{S^1}\mathcal{U}(p)$ in the form
$1\otimes(\gamma_1 \wedge \gamma_2)$. For
this, we compute a $\bbZ$-basis modulo $\{\pm 1\}$  of the
f.g.\ multiplicative abelian group $U_{S(p)}(K^+)$.
(Note that functions to perform this computation are implemented in PARI/GP.)
Once a basis is known, we use it to construct two random elements
$\gamma_1$, $\gamma_2$ in $U_{S(p)}(K^+)$. If $1\otimes (\gamma_1 \wedge
\gamma_2)$ does not lie in $e_{S^1}\cU(p)$
we replace, say, $\gamma_1$ by
$(|G|e_{S^1})\gamma_1$ so that it does.  Then  $1\otimes(\gamma_1 \wedge \gamma_2)$ will
generate $e_{S^1}\cU(p)$ if (and, in fact, only if) $\chi(R_{K^+/k,S(p)}((1\otimes\gamma_1)\wedge (1\otimes\gamma_1)))$ is non-zero for all characters $\chi
\in \hat{\bar{G}}$ such that $\ord_{s=0}L_{K^+/k, S}(s,\chi)=2$.
These conditions can be unconditionally tested using a good enough
approximation to $R_{K^+/k,S(p)}((1\otimes\gamma_1)\wedge (1\otimes\gamma_1))$, calculated as a group-ring determinant involving real logarithms of
(absolute values of) conjugates of $\gamma_1$ and $\gamma_2$. If they are not satisfied, we recommence with two new random elements
$\gamma_1$ and $\gamma_2$. (For our initial `random' choices of  $\gamma_1$ and $\gamma_2$ we actually took pairs of distinct elements of
the computed $\mathbb{Z}$-basis of $U_{S(p)}(K^+)$. In the few cases where this did not provide a generator,
%%  we then we looked at random `simple' elements.)
we then looked at pairs consisting of `simple' random linear combinations of these basis elements.)

We now know that the unique element $\eta_{K^+/k,S^1}$ -- if it exists -- will be equal to $A(1\otimes(\gamma_1 \wedge \gamma_2))$
for any $A \in \mathbb{Q}\bar{G}$ satisfying
\beq\label{eq: computation of  A} A \, R_{K^+/k,S^1}((1\otimes\gamma_1)\wedge (1\otimes\gamma_1))
=\Theta^{(d)}_{K^+/k,S^1}(0)
\eeq
(by~(\ref{eq: RS element})).
We compute an approximation of $\Theta^{(d)}_{K^+/k,S^1}(0)$ in $\mathbb{R}\bar{G}$ using its
expression in terms of Artin $L$-functions (see the beginning of
Section~\ref{sec: RS + pairing H}), once again using the methods of
\cite{D-T}, or \cite[Section 10.3]{Co}. Then we can find a solution $\tilde{A} \in
\mathbb{R}\bar{G}$ of  Equation~(\ref{eq:  computation of A}) to a high precision. Standard
methods allow us to compute an element $A_0 \in \mathbb{Q}\bar{G}$ very close to $\tilde{A}$
and with coefficients of small height. We then write $A_0$ as $\frac{1}{a} B_0$
where $a \in \mathbb{Z}_{>0}$ and $B_0$ is an element of  $\mathbb{Z}\bar{G}$, the g.c.d.\ of whose coefficients is prime to $a$.
Assuming that $A_0$ is in fact an {\em exact} solution of Equation~(\ref{eq:  computation of A}) (see below)
we now have the desired expression~(\ref{eq: what to expect 3}) with $\varepsilon_1 = \gamma_1$, $\varepsilon_2 = B_0 \, \gamma_2$.
However, we shall see in the next section that the computations of $H_{K/k,n}(\eta_{K^+/k,S^1},\theta)$ are {\em much}
easier if~(\ref{eq: what to expect 2}) holds. Thus if
$p$ divides $a$ we find a new generator $1\otimes (\gamma_1\wedge \gamma_2)$ and repeat the process.
We have justified above the expectation that~(\ref{eq: what to expect 2}) is possible whenever $p\ndiv|\bar{G}|$  and indeed the above process
terminated with such an expression in all our examples of this type. More surprisingly, perhaps, it also
terminated with a solution of~(\ref{eq: what to expect 2}) in all our examples with $p||\bar{G}|$.
Very similar behaviour was observed in~\cite{Solomon-Roblot1} (see also~\cite[Rem.\ 3.4]{zetap1}). Thus it seems, experimentally at least,  that
Rubin-Stark elements are `usually' better-behaved in this sense than the various conjectures predict, although no convincing
sharpening has yet been proposed along these lines.
\vertsp\\
\rem\ We need to convince ourselves that this is indeed the
Rubin-Stark element and not some \textit{ad hoc} element of
$e_{S^1}\mathcal{U}_{p}$ constructed simply to satisfy Equation~(\ref{eq: RS element}) to
the working precision. A first significant fact is that while we are
working with a large precision -- usually of $100$ digits -- the
coefficients of $A_0$ are of very small height. In almost all examples
numerators and denominators of the coefficients of $A_0$ are less than
$10$ in absolute value, the largest ones being in example \texttt{E5} where they have up to $6$
digits. However this is still considerably smaller than one would
expect if $\tilde{A}$ were a random element of $\mathbb{R}\bar{G}$. A
second and even more convincing way to reassure ourselves that we really have the Rubin-Stark element is as follows. Once
we have calculated an element $\hat{\eta}=\frac{1}{a}\otimes(\eps_1\wedge\eps_2)$, say, of $e_{S^1}\mathcal{U}(p)$ as a candidate for $\eta_{K^+/k,S^1}$, we
significantly increase the working precision, say from $100$ to $150$
digits. We then recompute $R_{K^+/k}(\hat{\eta})$ and
$\Theta^{(d)}_{K^+/k,S^1}(0)$ to the new precision and check whether they
still agree. If $\hat{\eta}$ were an {\em ad hoc} element,
constructed to satisfy Equation~(\ref{eq: RS element}) to a precision
of $100$ digits, then there would be no reason for it to satisfy it to $150$ digits. The fact that
it always did so, without readjustment, was, we felt, convincing
enough evidence to take $\eta_{K^+/k,S^1}=\hat{\eta}$.

\subsection{Computation of $H_{K/k,n}$ and Verification of the
  Conjecture}\label{subsec: compute H}
To complete the verification of the Congruence Conjecture, it suffices to check that $\sKk(\theta)$ lies in $\bbZ_p G^-$ and that the
two sides of~(\ref{eq: the congruence}) agree, for all $\theta$  in an appropriate subset of $\WUotwo$.
The determination of this subset, as well as the treatment of the L.H.S.,
is divided into three cases. In Case~$1$ (\ie\ examples \texttt{B6}, \texttt{C8}, \texttt{D7},
\texttt{D9} and \texttt{D11})  $\eta_{K^+/k,S^1}=0$ because $e_{S^1}=0$.
It then suffices to calculate an approximation to $\sKk(\theta)$ up to an element of $\pnpo\bbZ_p G^-$
for each $\theta$ in $W$, the initial $\bbZ_p G$-generating set for $\WUotwo$ constructed in Section~\ref{sec: gens of exterior square}.
These approximations are calculated using
Proposition~\ref{prop: convergence of sKk} with $M=n+1$ and the conjecture is verified if and only if each actually
lies in $\pnpo\bbZ_p G^-$. In the remainder of our examples, $\eta_{K^+/k,S^1}$ is {\em non-zero}
and the computation of the R.H.S.\ of~(\ref{eq: the congruence}) is usually lengthy.  In Case~$2$, $p\ndiv|G|$ and we explained at the beginning of
Subsection~\ref{sec: gens of exterior square} why it is sufficient to
check~(\ref{eq: the congruence}) for each $\theta$ in
the much smaller set
$W'$ generating modulo
$\WUotwo^+ +\left(\WUotwo\right)_{\rm tor}$ constructed there. This very construction
included the computation of an element of $\bbQ G^-$ approximating $\sKkS(\theta)$  up to an element of
$\pnpo\bbZ_p G^-$ (at least) for all $\theta\in W'$. In Case~$3$, $p||G|$ and we are unable to reduce $W$. So, once again we
use Proposition~\ref{prop: convergence of sKk} to
calculate an approximation to $\sKk(\theta)$ up to an element of $\pnpo\bbZ_p G^-$
for each $\theta$ in the full set $W$.

It remains to explain the computation of the R.H.S.\ of~(\ref{eq: the congruence}) in the second and third cases above, where
$\eta_{K^+/k}:=\eta_{K^+/k,S^1}\neq 0$.
The $\tau_i$ are realised as elements of $\Gal(F/\bbQ)$ and since $F$ contains $K$ and hence $\mu_{\pnpo}$,
the quantity $\kappa_n(\tau_1\tau_2)$ may be determined directly by calculating
$\xi_\pnpo^{\tau_1\tau_2}$. The computation of $H_{K/k,n}(\eta_{K^+/k},\theta)$
is greatly facilitated by the fact that Equation~(\ref{eq: what to expect 2}) -- hence also~(\ref{eq: what to expect}) -- holds in every
case, as already noted. Indeed, from diagram~(\ref{diag: H factors through H}) it follows that
\beq\label{eq: easier formula}
H_{K/k,n}(\eta_{K^+/k},\theta)=\bar{a}\inv
\cH_{K/k,n}(\eps_1\wedge\eps_2,\theta)\ \ \ \mbox{for all $\theta\in\WUotwo$}
\eeq
where $\bar{a}$ is the reduction of $a$ modulo $p^{n+1}$. Recall that every $\theta\in W$ is `global' by construction,
\ie\ of the form $\iota(v_1)\wedge\iota(v_2)$ for some $v_1,v_2\in K^\times$. Therefore,
using~(\ref{eq: easier formula}), the conditions satisfied by the
$\eps_i$ and the definitions of $\cH_{K/k,n}$
and $[\cdot,\cdot]_{K,n,G}$, it suffices to be able to calculate
$[\eps,\iota(v)]_{K,n}$ for any $v\in K^\times$ and $\eps\in U_{S(p)}(K^+)$. The next Proposition
shows how we did this. (The basic idea is well-known, see \eg~\cite[\S~II.7.5]{Gras}.)
Let $\fQ$ be any prime ideal of $\cO_K$. If $\fQ\nin S_p(K)$ then
reduction modulo $\fQ$ gives an injection $\mu_\pnpo(K)\rightarrow(\cO_K/\fQ)^\times$ so that
$p^{n+1}|(N\fQ-1)$ and the image is the subgroup of $(N\fQ-1)/\pnpo$-th powers in $(\cO_K/\fQ)^\times$.
Thus, for each such $\fQ$ there is a homomorphism
${\rm apr}_{\fQ,n}:(\cO_K/\fQ)^\times\rightarrow\ZpnpoZ$ (the {\em additive, $\pnpo$-th power residue symbol modulo $\fQ$})
uniquely defined by
$
\bar{\xi}_\pnpo^{{\rm apr}_{\fQ,n}(\bar{b})}=\bar{b}^{(N\fQ-1)/\pnpo}
$
for all $\bar{b}\in(\cO_K/\fQ)^\times$. For the small values of $\pnpo$ occurring here, ${\rm apr}_{\fQ,n}$ is quick to calculate directly
and we have:
\begin{prop}
If $\eps\in U_{S(p)}(K^+)$ and $v\in K^\times$, then
\[
[\eps,\iota(v)]_{K,n}=\sum_{\fQ\nin S_p(K)}\ord_\fQ(v){\rm apr}_{\fQ,n}(\bar{\eps})
\]
where $\fQ$ runs over the (finite) set of prime ideals $\fQ$ of $\cO_K$ not dividing
$p$ (and such that $\ord_\fQ(v)\neq 0$).
\end{prop}
\bPf\ Let $L^\eps$ be the Kummer extension $K(\eps^{1/\pnpo})$ and write $h_\eps$ for the isomorphism $\Gal(L^\eps/K)\rightarrow\ZpnpoZ$ given by
$h_\eps(g)=\Ind_n(g(\eps^{1/\pnpo})/\eps^{1/\pnpo})$. (Everything is independent of the choice of root $\eps^{1/\pnpo}$.)
For each prime ideal $\fQ$ of $\cO_K$  we choose a prime ideal $\tilde{\fQ}$ dividing $\fQ$ in $\cO_{L^\eps}$ and write $\rec_\fQ$
for the composite homomorphism
\[
K^\times\stackrel{\iota_\fQ}{\longrightarrow} K_\fQ^\times\longrightarrow D_{\tilde{\fQ}}(L^\eps/K)\hookrightarrow\Gal(L^\eps/K)
\]
where the second homomorphism is the local reciprocity map. (This is independent of the choice of $\tilde{\fQ}$.)
It follows easily from the definition
and alternating property of the local Hilbert symbol $(\cdot,\cdot)_{K_\fP,\pnpo}$ (see~\cite[Prop.~3.2]{NeukirchANT}) for $\fP\in S_p(K)$ that
\beq\label{eq: heps for fP}
h_\eps(\rec_\fP(v))=\Ind_n(\iota_\fP\inv(\iota_\fP(v),\iota_\fP(\eps)))=-\Ind_n(\iota_\fP\inv(\iota_\fP(\eps),\iota_\fP(v))_{K_\fP,\pnpo})\ \ \ \mbox{for all $\fP\in S_p(K)$}
\eeq
On the other hand, if $\fQ\nin S_p(K)$ then the extension $L^\eps/K$ is unramified  at $\fQ$ (since $\eps\in U_{S(p)}(K^+)$) and so
$\rec_\fQ(v)=\sigma_{\fQ,L^\eps/K}^{\ord_\fQ(v)}$ where $\sigma_{\fQ,L^\eps/K}$ denotes the Frobenius element. Since $\eps^{1/\pnpo}$ is a local unit at $\tilde{\fQ}$, the definition of
$\sigma_{\fQ,L^\eps/K}$ tells us that the image of the $\pnpo$th root of unity
$\sigma_{\fQ,L^\eps/K}(\eps^{1/\pnpo})/\eps^{1/\pnpo}$ in $(\cO_K/\fQ)^\times\subset(\cO_{L^\eps}/\tilde{\fQ})^\times$ is equal to that of
$\eps^{(N\fQ-1)/\pnpo}$. It follows that
\beq\label{eq: heps for fQ}
h_\eps(\rec_\fQ(v))=\ord_\fQ(v)h_\eps(\sigma_{\fQ,L^\eps/K})=\ord_\fQ(v){\rm apr}_{\fQ,n}(\bar{\eps})\ \ \ \mbox{for all $\fQ\nin S_p(K)$}
\eeq
In particular  $h_\eps(\rec_\fQ(v))$, and therefore $\rec_\fQ(v)$, is trivial for almost all $\fQ$.
Finally,
global class-field theory tells us that the product of $\rec_\fQ(v)$ over all prime ideals is equal to $1\in \Gal(L^\eps/K)$. (Since $K$ is totally complex the local reciprocity map is
trivial at archimedean places.) Using this and equations~(\ref{eq: def of square brackets}),  (\ref{eq: heps for fP}) and~(\ref{eq: heps for fQ}), we get
\[
[\eps,\iota(v)]_{K,n}=
%% \sum_{\fP\in S_p(K)}\Ind_n(\iota_\fP\inv(\iota_\fP(\eps),\iota_\fP(v))_{K_\fP,\pnpo})=
-\sum_{\fP\in S_p(K)}h_\eps(\rec_\fP(v))=
\sum_{\fQ\nin S_p(K)}h_\eps(\rec_\fQ(v))=\sum_{\fQ\nin S_p(K)}\ord_\fQ(v){\rm apr}_{\fQ,n}(\bar{\eps})
\]
as required. \ePf

\rem\label{rem: factor}\
In order to compute $H_{K/k,n}(\eta_{K^+/k}, \theta)$ for $\theta$ in $W$ (or $W'$) using the Proposition, one needs to compute
the prime-ideal factorisation of $(v_1)$ and $(v_2)$ in $K$ where $\theta
= \iota(v_1) \wedge \iota(v_2)$. Since, moreover, the $v_i$'s lie in $\cO_K$, the first step is to
factor the absolute norm of $v_i$, $i=1,2$. Unfortunately, the $v_i$'s  constructed by the method of Propositions~\ref{prop: gensofUone1}
and~\ref{prop: gensofUone2} tend to have very large norms which can be divisible by more than one large prime
number and hence virtually impossible to factor. We get around this problem by perturbing one or more of the $v_i$'s, \ie\ replacing
$v_i$ by $v'_i:=v_i + p^{n+2} x_i$ for a random element $x_i\in\cO_K$ for $i=1,2$.  Remark~\ref{rem: modifying v's} (with $l=n+1$)
implies that $\iota(v'_1) \wedge \iota(v'_2)\equiv\theta$ modulo $p^{n+1}\bigwedge^2_{\bbZ_p G}U^1(K_p)$ and so
$H_{K/k,n}(\eta_{K^+/k}, \theta)=H_{K/k,n}(\eta_{K^+/k}, \iota(v'_1) \wedge \iota(v'_2))$.
%% It follows firstly (by Nakayama's
%% Lemma) that such perturbations do not effect the generation of $\bigwedge^2_{\bbZ_p G}U^1(K_p)$ by  $W$ (or of
%% $\bigwedge^2_{\bbZ_p G}U^1(K_p)$ modulo $\bigwedge^2_{\bbZ_p G}U^1(K_p)^+ + \left(\bigwedge^2_{\bbZ_p G}U^1(K_p)\right)_{\rm tor}$ by  $W'$).
%% Secondly, we must have $\sKk(\iota(v'_1) \wedge \iota(v'_2))\equiv\sKk(\theta)$ modulo
%% $p^{n+1}\SKk$ and since by this point we know that
%% $\SKk\subset\bbZ_p G$, we conclude that the
Thus we may proceed as follows. We set some time
limit, say two minutes during which we try to factor the norm of
of each $v_i$. If we fail, we just perturb one or more of the $v_i$'s as above
and try again. Indeed, although $v'_i$'s  usually have norms of
about the same size as those of the $v_i$'s, it usually happens
that after several tries, we find norms that are (relatively) easy to factor,
allowing us to calculate $H_{K/k,n}(\eta_{K^+/k}, \iota(v'_1) \wedge \iota(v'_2))$
\ie\  $H_{K/k,n}(\eta_{K^+/k}, \theta)$.

% We use the following trick to get around this problem.
%  Recall that $\sKk(\theta)$ is computed using
%  a precision of $p^{-M}$ for a certain $M \geq n+1$ (see
%  Section~\ref{subsec: comp of sKk}), hence we can replace $v_i$ by $v_i
%  + p^M w_i$ for some element $w_i \in \mathcal{O}_K$, $i=1,2$, without
%  changing the value of $\sKk(\theta)$, and thus the value of the
%  L.H.S. of (\ref{eq: the congruence}). Therefore, we proceed as
%  follows: we set some time limit, say $5$ minutes, and for $i = 1, 2$,
%  we try to factor the norm of $v_i$ during that time. If we fail, we
%  just replace $v_i$ by $v_i + p^M w$ for some random element $w \in
%  \mathcal{O}_K$ and try again. Indeed, although the new element $v_i +
%  p^M w$ has usually a norm of similar size than that of $v_i$, it
%  usually happen that after several tries, we find a norm that is
%  easy to factor.
%
%
%  [*** in rephrasing the below one can say instead  that the $v_i$ can be replaced by $v_i
%  + p^M w_i$ for any elements element $w_i \in \mathcal{O}_K$, $i=1,2$, without
%  changing the value of $\sKk(\theta)$ etc... In this way we can ensure that they have small coefficients in terms of a computed $\bbZ$-basis
%  of $\cO_K$ (\eg\ of absolute value at most $p^M$: compute the coefficients of the original elements
%  in this basis, then take their least positive residues modulo $p^M$). Ask Xavier if this is OK. Perhaps it's even what he did. In
%  any case, it looks a bit better than `random modifications' which still
%  might given numbers with very large norms, as noted below...*****]
%
%
\section{Results of the Computations}\label{sec: results}

\subsection{An Example}

We illustrate the numerical computations with example \texttt{B1} (see
next subsection). We have $p = 3$ and $n = 0$, $k$ is the real
quadratic field $\mathbb{Q}(\sqrt{6})$ (thus $p$ ramifies in
$k/\mathbb{Q}$), $K^+$ is the ray-class field of $k$ of conductor
$4\mathfrak{p}$ where $\mathfrak{p}$ is the unique prime ideal of $k$
dividing $3$, and $K = K^+(\xi_3)$. The extension $K^+/k$ is of degree
$4$ with Galois group $\bar{G}$ isomorphic to $C_2^2$, and
the extension $K/k$ has degree $8$ and its Galois group $G$ is
isomorphic to $C_2^3$. In particular, $p$ does not divide
$|G|$.

The extension $K/\mathbb{Q}$ is a Galois extension, but is not
abelian, and we have $K = \mathbb{Q}(\nu)$ where $\nu$ is a root of the irreducible polynomial
\begin{multline*}
  X^{16} - 8 X^{15} + 48 X^{14} - 196 X^{13} + 642 X^{12} - 1668
  X^{11} + 3580 X^{10} - 6328 X^9 + 9297 X^8 \\
  - 11276 X^7 + 11224 X^6 - 9024 X^5 + 5736 X^4 - 2780 X^3 + 972 X^2 -
  220 X + 25
\end{multline*}

We find $\fp\cO_{K^+}=\fP_+^2$ (so that
%% The prime $p=3$ has the following decomposition in
%% $K^+$ and $K$ respectively
%% $$
%% 3\mathcal{O}_{K+} = \mathfrak{P}_+^4, \quad 3\mathcal{O}_K = \mathfrak{P}_1^4\mathfrak{P}_2^4
%% $$
$e_{\mathfrak{P}_+}(K^+/k) = f_{\mathfrak{P}_+}(K^+/k) = 2$) and
%% $e_{\mathfrak{P}_i}(K/K^+) = f_{\mathfrak{P}_i}(K/K^+) = 1$ for $i =
%% 1, 2$.
$\fP_+\cO_K=\fP\fP'$.
Finally, we have $S = S^1 = \{\infty_1, \infty_2, \mathfrak{p},
\mathfrak{q}_2\}$ where $2\mathcal{O}_k = \mathfrak{q}_2^2$.

Let $\sigma_1, \sigma_2, \sigma_3$ be three distinct $k$-automorphisms
of $K$ of order $2$ such that $G = \langle \sigma_1, \sigma_2, \sigma_3
\rangle$, with the convention that $\sigma_3 : \nu \mapsto 1-\nu$ is
the complex conjugation of the CM field $K$. Using the method of
Subsection~\ref{subsec: comp of sKk}, we find that
$$
a_{K/k}^- = \frac{1}{2^6 \, 3^2}\big(\sigma_3 - 1\big)
  \left(3 + 2\sqrt{3} + \sigma_1 + \sigma_1\sigma_2 + \big(3 -
   2\sqrt{3}\big) \sigma_2\right)
$$

With the notations of Subsection~\ref{subsec: generators of U1}, we
have $t = 1$, $h_1 = 2$ (with $\mathfrak{p}(1) = \mathfrak{p}$,
$\mathfrak{P}(1,1) = \mathfrak{P}$, $\mathfrak{P}(1,2) =
\mathfrak{P}'$) and $e_{\fP(1,j)}(K/\mathbb{Q})=4$ for $j=1,2$. Thus $l_1 = 7$, and Propositions~\ref{prop: gensofUone1}
and~\ref{prop: gensofUone2} enable us to construct $6$ elements such
that the set $W$ of wedge product of two of these generate $\WUotwo$
over $\mathbb{Z}_pG$. We now use the method (and the notations) of
Subsection~\ref{sec: gens of exterior square} to find a smaller
generating subset. Let $\chi_i$, $i=1,2,3$, be the character of $G$
defined by $\chi_i(\sigma_i) = -1$ and $\chi_i(\sigma_j) = 1$ for $j
\not= i$. It is easy to see that the set $R^- := \{\chi_3,
\chi_1\chi_3, \chi_2\chi_3, \chi_1\chi_2\chi_3\}$ is a system of
representatives of the orbits of the odd, irreducible characters of
$G$ under the action of
$\mathrm{Gal}(\bar{\mathbb{Q}}/\mathbb{Q})$. Thus, we have
$F_\chi=\bbQ$ for all $\chi$ in Equation~(\ref{eq: the iso alpha}),
and the equation gives
$$
\mathbb{Q}_pG^- \simeq \mathbb{Q}_p^4
$$
We compute that
$$
m(\chi_3, (p)) = m(\chi_2\chi_3, (p)) = m(\chi_1\chi_2\chi_3, (p)) = 0 \quad
\text{and} \quad m(\chi_1\chi_3, (p)) = 1
$$
and after several tries, we find a set $W$ such that
$\bigcap_{\chi \in R^-}W_{\min}(\chi,p)$ is non-empty and we take in
this set the element $\theta_0 = \iota(v_1) \wedge \iota(v_2)$ where
\begin{multline*}
  v_1 = \frac{1}{17095}\big(1058221 \nu^{15} - 7915486 \nu^{14} +
  46551510 \nu^{13} - 182313497 \nu^{12} + 579396826 \nu^{11} \\
  - 1444318673 \nu^{10} + 2976716004 \nu^9 - 5002660697 \nu^8 +
  6945207975 \nu^7 \\
  - 7851102425 \nu^6 + 7170233086 \nu^5 - 5155280183 \nu^4 +
  2822456537 \nu^3 \\
  - 1105885714 \nu^2 + 278328786 \nu - 33994775
  \big)
\end{multline*}
and
\begin{multline*}
  v_2 = \frac{1}{17095}\big(-383541 \nu^{15} + 2749923 \nu^{14} -
  16006808 \nu^{13} + 61029582 \nu^{12} - 190600453 \nu^{11} \\
  + 462662235 \nu^{10} - 930346920 \nu^9 + 1513004524 \nu^8 -
  2026236417 \nu^7 \\
  + 2184191092 \nu^6 - 1881836887 \nu^5 + 1247007651 \nu^4 - 609767073
  \nu^3 \\
  + 198580288 \nu^2 - 36118966 \nu + 1344335 \big)
\end{multline*}

By the result of Subsection~\ref{sec: gens of exterior square}, we
know that $\theta_0$ generates $\WUotwo$ over $\mathbb{Z}_p G$ modulo
$\WUotwo^++\left(\WUotwo\right)_{\rm tor}$ so to prove
$CC(K/k,S^1,p,n)$ it suffices to establish~(\ref{eq: the congruence})
with $\theta=\theta_0$. Note that the L.H.S.\ of ~(\ref{eq: the
  congruence}) has already been computed. For the R.H.S.,\ the field
$K^+$ is generated over $\mathbb{Q}$ by a root $\lambda$ of the
irreducible polynomial
$$
P_\lambda(X) = X^8 - 4 X^7 - 4 X^6 + 20 X^5 + 4 X^4 - 20 X^3 - 4 X^2 +
4 X + 1
$$
Using the methods described in Subsection~\ref{subsec: compute
  rs-element}, we find that the Rubin-Stark element is given by
$$
\eta_{K^+/k} = \frac{1}{16} (\eps_1 \wedge \eps_2)
$$
where
$$
\eps_1 = \frac{1}{5}\big( 6\lambda^7 - 22\lambda^6 - 33\lambda^5 +
119\lambda^4 + 52\lambda^3 - 121\lambda^2 - 31\lambda + 7 \big)
$$
and
\begin{multline*}
  \eps_2 = \frac{1}{25}\big(-102282\lambda^7 + 463929\lambda^6 +
  152556\lambda^5 - 2073598\lambda^4 \\
  + 604836\lambda^3 + 1722767\lambda^2 - 413178\lambda - 221449 \big).
\end{multline*}
Note that $\eps_1$ and $\eps_2$ lie in $\cO^\times_{K^+}$, and not just
in $U_{S(p)}(K^+)$.

We now compute by Subsection~\ref{subsec: compute H}\footnote{As
  mentioned in Remark~\ref{rem: factor}, one needs to factor the norm
  of $v_1$ and $v_2$ to do this computation, but since these are of
  about $17$ digits, it is easy in this case.}
$$
H_{K/k,0}(\eta_{K^+/k}, \theta_0) = (\sigma_3 - \bar{1})(\sigma_1 -
\sigma_2 - \bar{1}) \in (\bbZ/3\bbZ)G.
$$
Finally, we can check that
$$
\sKk(\theta_0) \equiv H_{K/k,0}(\eta_{K^+/k}, \theta_0) \pmod{3}
$$
and therefore $CC(K/k,S^1,p,n)$ is satisfied (since we compute that
our choice of $\tau_1$ and $\tau_2$ implies that
$\kappa_0(\tau_1\tau_2) \equiv 1 \pmod{3}$).

\subsection{Tables}

We have numerically verified that the conjecture $CC$ is satisfied in
$48$ examples. These examples are divided into $4$
types\footnote{A fifth, type \texttt{A}, for which the extension
$K/\mathbb{Q}$ is abelian, was used for testing purposes only. It is not included
because the $CC$ then follows from~\cite[Thm.~5]{laser}. (Hypothesis~4, {\em ibid.}
holds since $p\ndiv 2=[k:\bbQ]$).} of differing significance.
\begin{itemize}
\item $12$ examples of type \texttt{B}: $p = 3$, $5$ or $7$, $n = 0$,
  $p$ does not divide $|G|$, $K/\mathbb{Q}$ is Galois but not abelian;
\item $16$ examples of type \texttt{C}: $p = 3$, $5$ or $7$, $n = 0$,
  $p$ does not divide $|G|$, $K/\mathbb{Q}$ is non-Galois;
\item $14$ examples of type \texttt{D}: $p = 3$ or $5$, $n = 0$, $p$
  divides $|G|$, $K/\mathbb{Q}$ non-Galois (resp. Galois but not
  abelian) if $p = 3$ (resp. $p= 5$);
\item $6$ examples of type \texttt{E}: $p = 3$, $n = 1$, $p$
  necessarily divides $|G|$, $K/\mathbb{Q}$ not abelian but possibly
  Galois.
\end{itemize}

The examples are summarized in four tables, with one table for each
type. The columns of the tables have the following meaning:
\begin{itemize}
\item the number of the example,
\item the value of $p$ (it is either $3$, $5$ or $7$),
\item the discriminant $d_k$ of the real quadratic base field $k$
  (thus $k = \mathbb{Q}(\sqrt{d_k})$),
\item `$R$', `$S$' or `$I$' according to whether $p$ is ramified,
  split or inert in $k$,
\item the conductor $\mathfrak{f}(K)$ of $K/k$ with the following
  notations: $\fp = \fp(1)$, and $\fp' = \fp(2)$ if $p$ is split in
  $k$, $\fq_q$ is a prime ideal of $k$ above a primer number $q$, and
  $\fq_q'$ is the other prime ideal of $k$ above $q$ if $q$ is split
  in $k$,\footnote{Note that the examples \texttt{B8} and \texttt{B9}
    differ only by the prime ideal in $k$ above $7$ dividing the
    conductor.}
\item the structure of the Galois group $G$ as a product of
  cyclic groups,
\item the structure of the Galois group $\bar{G}$ as a product of
  cyclic groups,
\item the minimal polynomial $P_\lambda$ of a generating element
  $\lambda$ of $K^+$ over $\mathbb{Q}$, so $K^+ = \mathbb{Q}(\lambda)$
  and $K = \mathbb{Q}(\lambda, \xi_p)$,
\item the decomposition in $K/k$ of the primes ideals above $p$
  given as $(e_{\fP(i,1)}(K/k), f_{\fP(i,1)}(K/k), h_i)$ for
  $i=1,\dots,t$ (see Subsection~\ref{subsec: generators of U1} for the
  notations),
\item the cardinality of $S^1$,
\item the value of $a$ (see Equation~\eqref{eq: what to expect 2}),
\item the nature of the Rubin-Stark element: a `$0$' means that it is
  trivial, a `$U$' means that we found a representation as
  in~\eqref{eq: what to expect 2} with
  $\varepsilon_i\in\cO_{K^+}^\times$ for $i=1,2$. Recall from
  Subsection~3.5 that this is to be expected in examples where $p\ndiv
  |G|$ and $|S^1|\geq 4$. Interestingly, it turned out to be possible
  in most of our other examples as well. In the remainder, indicated
  by a `$p$' in this column, we were only able to satisfy~\eqref{eq:
    what to expect 2} with $\varepsilon_1\in\cO_{K^+}^\times$ and
  $\varepsilon_2\in U_{S(p)}(K^+)$.
\end{itemize}

\begin{sidewaystable}
\centering
{\fontsize{8pt}{9pt}\selectfont
\begin{tabular}{|c|c|c|c|c|c|c|c|c|c|c|c|}
  \hline
  \onecolc{12}{\footnotesize \textrm{Table 1: Examples of type
      \texttt{B} ($n = 0$, $p \nmid |G|$, $K/\mathbb{Q}$ Galois but
      not abelian)
      \phv{12}{6}}} \\
  \hline
  \# & \phv{11}{5}
  $p$ & $d_k$  & $p$ in $k$ & $\ff(K)$ & $G$ & $\bar{G}$ & $P_\lambda(X)$
  & $p$ in $K/k$ & $|S^1|$ & $a$ & $\eta_{K^+/k}$ \\
  \hline
  1  \phv{11}{5}
  & $3$
  & $24$
  & $R$
  & $4\fp$
  & $C_2^3$
  & $C_2^2$
  & \parbox{8cm}{
    \vskip-12pt \begin{multline*}
      X^8 - 4 X^7 - 4 X^6 + 20 X^5 + 4 X^4 - 20 X^3 - 4 X^2 + 4 X + 1
    \end{multline*}\ \vskip-12pt
    }
   & $(2,2,2)$
   & $4$
   & $16$
   & $U$
   \\
\hline
2  \phv{11}{5}
   & $3$
   & $28$
   & $S$
   & $3\fq_2^5$
   & $C_2^3$
   & $C_2^2$
   & \parbox{8cm}{
    \vskip-12pt \begin{multline*}
      X^8 - 4 X^7 - 8 X^6 + 24 X^5 + 30 X^4 - 16 X^3 - 20 X^2 + 2
    \end{multline*}\ \vskip-12pt
    }
   & $(2,2,2)(2,2,2)$
   & $5$
   & $2$
   & $U$
   \\
\hline
3  \phv{11}{5}
   & $3$
   & $29$
   & $I$
   & $15\cO_k$
   & $C_2^3$
   & $C_2^2$
   & \parbox{8cm}{
    \vskip-12pt \begin{multline*}
      X^8 - 2 X^7 - 12 X^6 + 26 X^5 + 17 X^4 - 36 X^3 - 5 X^2 + 11 X -
      1
    \end{multline*}\ \vskip-12pt
    }
   & $(2,2,2)$
   & $5$
   & $1$
   & $U$
   \\
\hline
4  \phv{11}{5}
   & $3$
   & $33$
   & $R$
   & $4\fp$
   & $C_2^3$
   & $C_2^2$
   & \parbox{5cm}{
     \vskip-12pt \begin{multline*}
      X^8 - 11 X^6 + 24 X^4 - 11 X^2 + 1
    \end{multline*}\ \vskip-12pt
    }
   & $(2,2,2)$
   & $5$
   & $4$
   & $U$
    \\
\hline
5  \phv{11}{5}
   & $3$
   & $40$
   & $S$
   & $3\fq_2^3\fq_5$
   & $C_8 \times C_2$
   & $C_8$
   & \parbox{11cm}{
    \vskip-12pt \begin{multline*}
      X^{16} - 4 X^{15} - 16 X^{14} + 76 X^{13} + 46 X^{12} - 392
      X^{11} \\ - 24 X^{10} + 928 X^9 - 23 X^8 - 1128 X^7 - 44 X^6 +
      672 X^5 \\ + 96 X^4 - 156 X^3 - 36 X^2 + 4 X + 1
    \end{multline*}\ \vskip-12pt
    }
   & $(2,8,1)(2,8,1)$
   & $6$
   & $2$
   & $U$
    \\
\hline
6  \phv{11}{5}
   & $3$
   & $41$
   & $I$
   & $15\cO_k$
   & $C_4 \times C_2$
   & $C_4$
   & \parbox{5cm}{
    \vskip-12pt \begin{multline*}
      X^8 - 19 X^6 + 41 X^4 - 19 X^2 + 1
    \end{multline*}\ \vskip-12pt
    }
   & $(2,1,4)$
   & $5$
   & $1$
   & $0$
   \\
\hline
7  \phv{11}{5}
   & $3$
   & $44$
   & $I$
   & $3\cO_k$
   & $C_4 \times C_2$
   & $C_4$
   & \parbox{5cm}{
    \vskip-12pt \begin{multline*}
      X^8 - 11 X^6 + 24 X^4 - 11 X^2 + 1
    \end{multline*}\ \vskip-12pt
    }
   & $(4,1,2)$
   & $3$
   & $80$
   & $p$
   \\
\hline
8  \phv{11}{5}
   & $3$
   & $505$
   & $S$
   & $3\cO_k$
   & $C_8 \times C_2$
   & $C_8$
   & \parbox{11cm}{
    \vskip-12pt \begin{multline*}
      X^{16} - 3 X^{15} - 50 X^{14} + 157 X^{13} + 800 X^{12} - 3014
      X^{11} - 4242 X^{10} \\ + 25193 X^9 - 6314 X^8 - 82099 X^7 +
      99216 X^6 + 38525 X^5 \\ - 125349 X^4 + 50387 X^3 + 19768 X^2 -
      14926 X + 2029
      \end{multline*}\ \vskip-12pt
    }
   & $(2,8,1)(2,8,1)$
   & $4$
   & $16$
   & $U$
   \\
\hline
9  \phv{11}{5}
   & $5$
   & $29$
   & $S$
   & $5\cO_k$
   & $C_4 \times C_2$
   & $C_2 \times C_2$
   & \parbox{8cm}{
    \vskip-12pt \begin{multline*}
      X^8 - 2 X^7 - 12 X^6 + 26 X^5 + 17 X^4 - 36 X^3 - 5 X^2 + 11 X -
      1
    \end{multline*}\ \vskip-12pt
    }
   & $(4,2,1)(4,2,1)$
   & $4$
   & $1$
   & $U$
   \\
\hline
10 \phv{11}{5}
   & $5$
   & $41$
   & $S$
   & $5\cO_k$
   & $C_4 \times C_2$
   & $C_4$
   & \parbox{5cm}{
    \vskip-12pt \begin{multline*}
      X^8 - 19 X^6 + 41 X^4 - 19 X^2 + 1
    \end{multline*}\ \vskip-12pt
    }
   & $(4,1,2)(4,1,2)$
   & $4$
   & $4$
   & $U$
   \\
\hline
11 \phv{11}{5}
   & $5$
   & $44$
   & $S$
   & $15\cO_k$
   & $C_4^2$
   & $C_4 \times C_2$
   & \parbox{11cm}{
    \vskip-12pt \begin{multline*}
      X^{16} - 33 X^{14} + 289 X^{12} - 990 X^{10} + 1470 X^8 - 990
      X^6 + 289 X^4 - 33 X^2 + 1
    \end{multline*}\ \vskip-12pt
    }
   & $(4,4,1)(4,4,1)$
   & $5$
   & $1$
   & $U$
   \\
\hline
12 \phv{11}{5}
   & $7$
   & $29$
   & $S$
   & $35\cO_k$
   & $C_6 \times C_2^2$
   & $C_6 \times C_2$
   & \parbox{11cm}{
    \vskip-12pt \begin{multline*}
      X^{24} - 8 X^{23} - 22 X^{22} + 308 X^{21} - 94 X^{20} - 4452
      X^{19} + 5808 X^{18} + 31789 X^{17} \\ - 59220 X^{16} - 122740
      X^{15} + 288660 X^{14} + 258142 X^{13} - 781464 X^{12} \\ -
      270957 X^{11} + 1221211 X^{10} + 92820 X^9 - 1092490 X^8 +
      34676 X^7 \\ + 537022 X^6 - 9659 X^5 - 133103 X^4 - 11639 X^3 +
      12521 X^2 + 2860 X + 169
    \end{multline*}\ \vskip-12pt
    }
   & $(6,2,2)(6,2,2)$
   & $6$
   & $3$
   & $U$
   \\
\hline
\end{tabular}
}
\end{sidewaystable}

\begin{sidewaystable}
\centering
{\fontsize{8pt}{9pt}\selectfont
\begin{tabular}{|c|c|c|c|c|c|c|c|c|c|c|c|}
  \hline
  \onecolc{12}{\footnotesize \textrm{Table 2: Examples of type
      \texttt{C} ($n = 0$, $p \nmid |G|$, $K/\mathbb{Q}$ not Galois)
      \phv{12}{6}}} \\
  \hline
  \# & \phv{11}{5}
  $p$ & $d_k$  & $p$ in $k$ & $\ff(K)$ & $G$ & $\bar{G}$ &
  $P_\lambda(X)$ & $p$ in $K/k$ & $|S^1|$ & $a$ & $\eta_{K^+/k}$ \\
  \hline
  1  \phv{11}{5}
  & $3$
  & $5$
  & $I$
  & $3\fq_{29}$
  & $C_2^2$
  & $C_2$
  & \parbox{4cm}{
    \vskip-12pt \begin{multline*}
      X^4 - X^3 - 3X^2 + X + 1
    \end{multline*}\ \vskip-12pt
    }
   & $(2,2,1)$
   & $4$
   & $1$
   & $U$
   \\
\hline
2  \phv{11}{5}
   & $3$
   & $8$
   & $I$
   & $3\fq_{41}$
   & $C_2^2$
   & $C_2$
   & \parbox{4cm}{
    \vskip-12pt \begin{multline*}
      X^4 - 2 X^3 - 3 X^2 + 2 X + 1
    \end{multline*}\ \vskip-12pt
    }
   & $(2,2,1)$
   & $4$
   & $1$
   & $U$
   \\
\hline
3  \phv{11}{5}
   & $3$
   & $12$
   & $R$
   & $\fp\fq_{11}$
   & $C_2^2$
   & $C_2$
   & \parbox{4cm}{
    \vskip-12pt \begin{multline*}
      X^4 - 2 X^3 - 3 X^2 + 4 X + 1
    \end{multline*}\ \vskip-12pt
    }
   & $(2,2,1)$
   & $4$
   & $2$
   & $U$
   \\
\hline
4  \phv{11}{5}
   & $3$
   & $13$
   & $S$
   & $12\cO_k$
   & $C_2^2$
   & $C_2$
   & \parbox{2.5cm}{
     \vskip-12pt \begin{multline*}
       X^4 - 5 X^2 + 3
    \end{multline*}\ \vskip-12pt
    }
   & $(2,1,2)(2,2,1)$
   & $5$
   & $1$
   & $U$
    \\
\hline
5  \phv{11}{5}
   & $3$
   & $17$
   & $I$
   & $3\fq_2^3\fq_2'^2$
   & $C_2^2$
   & $C_2$
   & \parbox{2.5cm}{
    \vskip-12pt \begin{multline*}
      X^4 - 5 X^2 + 2
    \end{multline*}\ \vskip-12pt
    }
   & $(2,2,1)$
   & $5$
   & $1$
   & $U$
    \\
\hline
6  \phv{11}{5}
   & $3$
   & $21$
   & $R$
   & $\fq_{37}$
   & $C_2^2$
   & $C_2$
   & \parbox{4cm}{
    \vskip-12pt \begin{multline*}
      X^4 - 2 X^3 - 4 X^2 + 5 X + 1
    \end{multline*}\ \vskip-12pt
    }
   & $(1,2,2)$
   & $4$
   & $1$
   & $U$
   \\
\hline
7  \phv{11}{5}
   & $3$
   & $24$
   & $R$
   & $4\fp$
   & $C_2^2$
   & $C_2$
   & \parbox{2.5cm}{
    \vskip-12pt \begin{multline*}
      X^4 - 6 X^2 + 3
    \end{multline*}\ \vskip-12pt
    }
   & $(2,1,2)$
   & $4$
   & $2$
   & $U$
   \\
\hline
8  \phv{11}{5}
   & $3$
   & $28$
   & $S$
   & $3\fq_2^5$
   & $C_2^2$
   & $C_2$
   & \parbox{2.5cm}{
    \vskip-12pt \begin{multline*}
      X^4 - 6 X^2 + 2
      \end{multline*}\ \vskip-12pt
    }
   & $(2,2,1)(2,1,2)$
   & $5$
   & $1$
   & $0$
   \\
\hline
9  \phv{11}{5}
   & $3$
   & $29$
   & $I$
   & $3\fq_5$
   & $C_2^2$
   & $C_2$
   & \parbox{4cm}{
    \vskip-12pt \begin{multline*}
      X^4 - X^3 - 5 X^2 - X + 1
    \end{multline*}\ \vskip-12pt
    }
   & $(2,2,1)$
   & $4$
   & $2$
   & $U$
   \\
\hline
10 \phv{11}{5}
   & $5$
   & $5$
   & $R$
   & $\fp\fq_{29}$
   & $C_2^2$
   & $C_2$
   & \parbox{4cm}{
    \vskip-12pt \begin{multline*}
      X^4 - X^3 - 3 X^2 + X + 1
    \end{multline*}\ \vskip-12pt
    }
   & $(2,2,1)$
   & $4$
   & $1$
   & $U$
   \\
\hline
11 \phv{11}{5}
   & $5$
   & $8$
   & $I$
   & $5\fq_{41}$
   & $C_4 \times C_2$
   & $C_2 \times C_2$
   & \parbox{8cm}{
    \vskip-12pt \begin{multline*}
      X^8 - 14 X^6 - 4 X^5 + 43 X^4 + 8 X^3 - 34 X^2 - 8 X + 4
    \end{multline*}\ \vskip-12pt
    }
   & $(4,1,2)$
   & $4$
   & $1$
   & $U$
   \\
\hline
12 \phv{11}{5}
   & $5$
   & $12$
   & $I$
   & $5\fq_3\fq_{11}$
   & $C_4 \times C_2$
   & $C_2 \times C_2$
   & \parbox{5cm}{
    \vskip-12pt \begin{multline*}
      X^8 - 14 X^6 + 45 X^4 - 32 X^2 + 4
    \end{multline*}\ \vskip-12pt
    }
   & $(4,2,1)$
   & $5$
   & $1$
   & $U$
   \\
\hline
13 \phv{11}{5}
   & $5$
   & $17$
   & $I$
   & $5\fq_2^3\fq_2'^2$
   & $C_4 \times C_2$
   & $C_2 \times C_2$
   & \parbox{5cm}{
    \vskip-12pt \begin{multline*}
      X^8 - 15 X^6 + 39 X^4 - 30 X^2 + 4
    \end{multline*}\ \vskip-12pt }
   & $(4,2,1)$
   & $5$
   & $1$
   & $U$
   \\
\hline
14 \phv{11}{5}
   & $7$
   & $5$
   & $I$
   & $\fp\fq_{29}$
   & $C_6 \times C_2$
   & $C_6$
   & \parbox{11cm}{
    \vskip-12pt \begin{multline*}
      X^{12} - X^{11} - 17 X^{10} + 8 X^9 + 79 X^8 - 32 X^7 - 126 X^6
      + 37 X^5 \\
      + 81 X^4 - 15 X^3 - 19 X^2 + 2 X + 1
    \end{multline*}\ \vskip-12pt
    }
   & $(6,1,2)$
   & $4$
   & $1$
   & $U$
   \\
\hline
15 \phv{11}{5}
   & $7$
   & $8$
   & $S$
   & $\fp\fq_{41}$
   & $C_6 \times C_2$
   & $C_6$
   & \parbox{11cm}{
    \vskip-12pt \begin{multline*}
      X^{12} - 2 X^{11} - 21 X^{10} + 30 X^9 + 142 X^8 - 146 X^7 - 383
      X^6 + 276 X^5 \\
      + 385 X^4 - 214 X^3 - 124 X^2 + 56 X - 1
    \end{multline*}\ \vskip-12pt
    }
   & $(6,2,1)(6,2,1)$
   & $5$
   & $1$
   & $U$
   \\
\hline
16 \phv{11}{5}
   & $7$
   & $28$
   & $R$
   & $\fp\fq_2^5$
   & $C_6 \times C_2$
   & $C_6$
   & \parbox{9cm}{
    \vskip-12pt \begin{multline*}
      X^{12} - 16 X^{10} + 88 X^8 - 204 X^6 + 212 X^4 - 88 X^2 + 8
    \end{multline*}\ \vskip-12pt
    }
   & $(3,2,2)$
   & $4$
   & $6$
   & $U$
   \\
\hline
\end{tabular}
}
\end{sidewaystable}

\begin{sidewaystable}
\centering
{\fontsize{8pt}{9pt}\selectfont
\begin{tabular}{|c|c|c|c|c|c|c|c|c|c|c|c|}
  \hline
  \onecolc{12}{\footnotesize \textrm{Table 3: Examples of type
      \texttt{D} ($n = 0$, $p \mid |G|$, $K/\mathbb{Q}$ not Galois if
      $p = 3$ and $K/\mathbb{Q}$ Galois but not abelian if $p = 5$)
      \phv{12}{6}}} \\
  \hline
  \# & \phv{11}{5}
  $p$ & $d_k$  & $p$ in $k$ & $\ff(K)$ & $G$ & $\bar{G}$ &
  $P_\lambda(X)$ & $p$ in $K/k$ & $|S^1|$ & $a$ & $\eta_{K^+/k}$ \\
  \hline
  1  \phv{11}{5}
  & $3$
  & $5$
  & $I$
  & $6\fq_{19}$
  & $C_6$
  & $C_3$
  & \parbox{6cm}{
    \vskip-12pt \begin{multline*}
      X^6 - X^5 - 6 X^4 + 7 X^3 + 4 X^2 - 5 X + 1
    \end{multline*}\ \vskip-12pt
    }
   & $(2,3,1)$
   & $5$
   & $1$
   & $U$
   \\
\hline
2  \phv{11}{5}
   & $3$
   & $8$
   & $I$
   & $3\fq_{79}$
   & $C_6$
   & $C_3$
   & \parbox{6cm}{
    \vskip-12pt \begin{multline*}
      X^6 - 2 X^5 - 5 X^4 + 10 X^3 - 4 X + 1
    \end{multline*}\ \vskip-12pt
    }
   & $(2,3,1)$
   & $4$
   & $1$
   & $U$
   \\
\hline
3  \phv{11}{5}
   & $3$
   & $12$
   & $R$
   & $15\cO_k$
   & $C_6$
   & $C_3$
   & \parbox{5cm}{
    \vskip-12pt \begin{multline*}
      X^6 - 21 X^4 - 10 X^3 + 42 X^2 - 8
    \end{multline*}\ \vskip-12pt
    }
   & $(3,2,1)$
   & $4$
   & $1$
   & $p$
   \\
\hline
4  \phv{11}{5}
   & $3$
   & $29$
   & $I$
   & $6\fq_7$
   & $C_6$
   & $C_3$
   & \parbox{7cm}{
     \vskip-12pt \begin{multline*}
       X^6 - 3 X^5 - 8 X^4 + 19 X^3 + 24 X^2 - 29 X - 29
    \end{multline*}\ \vskip-12pt
    }
   & $(2,3,1)$
   & $5$
   & $1$
   & $U$
    \\
\hline
5  \phv{11}{5}
   & $3$
   & $29$
   & $I$
   & $6\fq_{13}$
   & $C_6$
   & $C_3$
   & \parbox{7cm}{
    \vskip-12pt \begin{multline*}
      X^6 - X^5 - 14 X^4 - 13 X^3 + 6 X^2 + 7 X + 1
    \end{multline*}\ \vskip-12pt
    }
   & $(2,3,1)$
   & $5$
   & $1$
   & $U$
    \\
\hline
6  \phv{11}{5}
   & $3$
   & $29$
   & $I$
   & $3\fq_7\fq_{13}$
   & $C_6$
   & $C_3$
   & \parbox{7cm}{
    \vskip-12pt \begin{multline*}
      X^6 - 2 X^5 - 15 X^4 + 6 X^3 + 45 X^2 + 22 X - 4
    \end{multline*}\ \vskip-12pt
    }
   & $(2,3,2)$
   & $5$
   & $1$
   & $U$
   \\
\hline
7  \phv{11}{5}
   & $3$
   & $29$
   & $I$
   & $3\fq_7\fq_{13}$
   & $C_6$
   & $C_3$
   & \parbox{7cm}{
    \vskip-12pt \begin{multline*}
      X^6 - 2 X^5 - 16 X^4 + 6 X^3 + 76 X^2 + 79 X + 23
    \end{multline*}\ \vskip-12pt
    }
   & $(2,1,3)$
   & $5$
   & $1$
   & $0$
   \\
\hline
8  \phv{11}{5}
   & $3$
   & $37$
   & $S$
   & $\fp\fp'^2\fq_7$
   & $C_6$
   & $C_3$
   & \parbox{7cm}{
    \vskip-12pt \begin{multline*}
      X^6 - X^5 - 16 X^4 + 9 X^3 + 65 X^2 - 10 X - 4
    \end{multline*}\ \vskip-12pt
    }
   & $(2,3,1)(6,1,1)$
   & $5$
   & $1$
   & $U$
   \\
\hline
9  \phv{11}{5}
   & $3$
   & $37$
   & $S$
   & $\fp\fp'^2\fq_7'$
   & $C_6$
   & $C_3$
   & \parbox{7cm}{
    \vskip-12pt \begin{multline*}
      X^6 - X^5 - 14 X^4 + 20 X^3 + 16 X^2 - 16 X - 9
    \end{multline*}\ \vskip-12pt
    }
   & $(2,1,3)(6,1,1)$
   & $5$
   & $1$
   & $0$
   \\
\hline
10 \phv{11}{5}
   & $3$
   & $37$
   & $S$
   & $\fp\fq_{73}$
   & $C_6$
   & $C_3$
   & \parbox{7cm}{
    \vskip-12pt \begin{multline*}
      X^6 - 2 X^5 - 15 X^4 + 13 X^3 + 30 X^2 - 13 X - 7
    \end{multline*}\ \vskip-12pt
    }
   & $(2,3,1)(2,3,1)$
   & $5$
   & $1$
   & $U$
   \\
\hline
11 \phv{11}{5}
   & $3$
   & $40$
   & $S$
   & $\fp\fq_{37}$
   & $C_6$
   & $C_3$
   & \parbox{7cm}{
    \vskip-12pt \begin{multline*}
      X^6 - 2 X^5 - 9 X^4 + 18 X^3 + 12 X^2 - 20 X - 9
    \end{multline*}\ \vskip-12pt
    }
   & $(2,3,1)(2,1,3)$
   & $5$
   & $1$
   & $0$
   \\
\hline
12 \phv{11}{5}
   & $3$
   & $40$
   & $S$
   & $3\fq_{67}$
   & $C_6$
   & $C_3$
   & \parbox{7cm}{
    \vskip-12pt \begin{multline*}
      X^6 - 2 X^5 - 19 X^4 + 10 X^3 + 100 X^2 + 100 X + 25
    \end{multline*}\ \vskip-12pt
    }
   & $(2,3,1)(2,3,1)$
   & $5$
   & $1$
   & $U$
   \\
\hline
13 \phv{11}{5}
   & $3$
   & $44$
   & $I$
   & $3\fq_{19}$
   & $C_6$
   & $C_3$
   & \parbox{7cm}{
    \vskip-12pt \begin{multline*}
      X^6 - 2 X^5 - 10 X^4 + 14 X^3 + 19 X^2 - 22 X + 5
    \end{multline*}\ \vskip-12pt }
   & $(2,3,1)$
   & $4$
   & $1$
   & $U$
   \\
\hline
14 \phv{11}{5}
   & $5$
   & $5$
   & $R$
   & $55\cO_k$
   & $C_{10}$
   & $C_5$
   & \parbox{8cm}{
    \vskip-12pt \begin{multline*}
      X^{10} - 45 X^8 + 700 X^6 - 4265 X^4 + 7725 X^2 - 980
    \end{multline*}\ \vskip-12pt
    }
   & $(10,1,1)$
   & $5$
   & $1$
   & $U$
   \\
\hline
\end{tabular}
}
\end{sidewaystable}

\begin{sidewaystable}
\centering
{\fontsize{8pt}{9pt}\selectfont
\begin{tabular}{|c|c|c|c|c|c|c|c|c|c|c|c|}
  \hline
  \onecolc{12}{\footnotesize \textrm{Table 4: Examples of type
      \texttt{E} ($n = 1$ so $p \mid |G|$, $K/\mathbb{Q}$ not Galois
      or Galois but not abelian)
      \phv{12}{6}}} \\
  \hline
  \# & \phv{11}{5}
  $p$ & $d_k$  & $p$ in $k$ & $\ff(K)$ & $G$ & $\bar{G}$ &
  $P_\lambda(X)$ & $p$ in $K/k$ & $|S^1|$ & $a$ & $\eta_{K^+/k}$ \\
  \hline
  1  \phv{11}{5}
  & $3$
  & $13$
  & $S$
  & $90\cO_k$
  & $C_6 \times C_3$
  & $C_3^2$
  & \parbox{10cm}{
    \vskip-12pt \begin{multline*}
      X^{18} - 9 X^{17} - 3 X^{16} + 222 X^{15} - 387 X^{14} - 1701
      X^{13} + 4637 X^{12} \\ + 4659 X^{11} - 19920 X^{10} - 2160 X^9
      + 37413 X^8 - 5949 X^7 \\ - 32755 X^6 + 5439 X^5 + 12117 X^4 -
      441 X^3 - 813 X^2 - 60 X - 1
    \end{multline*}\ \vskip-12pt
    }
   & $(6,3,1)(6,3,1)$
   & $6$
   & $7$
   & $U$
   \\
\hline
2  \phv{11}{5}
   & $3$
   & $29$
   & $I$
   & $9\fq_5$
   & $C_6 \times C_2$
   & $C_6$
   & \parbox{10cm}{
    \vskip-12pt \begin{multline*}
      X^{12} - 3 X^{11} - 24 X^{10} + 52 X^9 + 195 X^8 - 306 X^7 - 680
      X^6 + 744 X^5 \\ + 999 X^4 - 727 X^3 - 516 X^2 + 213 X - 1
    \end{multline*}\ \vskip-12pt
    }
   & $(6,2,1)$
   & $4$
   & $1$
   & $U$
   \\
\hline
3  \phv{11}{5}
   & $3$
   & $37$
   & $S$
   & $18\cO_k$
   & $C_6 \times C_3$
   & $C_3^2$
   & \parbox{10cm}{
    \vskip-12pt \begin{multline*}
      X^{18} - 9 X^{17} + 3 X^{16} + 174 X^{15} - 357 X^{14} - 1083
      X^{13} + 3463 X^{12} \\ + 2001 X^{11} - 13218 X^{10} + 3150 X^9
      + 22479 X^8 - 14883 X^7 \\ - 15063 X^6 + 16155 X^5 + 741 X^4 -
      5073 X^3 + 1557 X^2 - 37
    \end{multline*}\ \vskip-12pt
    }
   & $(6,3,1)(6,3,1)$
   & $5$
   & $2$
   & $U$
   \\
\hline
4  \phv{11}{5}
   & $3$
   & $37$
   & $S$
   & $18\cO_k$
   & $C_6 \times C_3$
   & $C_3^2$
   & \parbox{10cm}{
     \vskip-12pt \begin{multline*}
       X^{18} - 9 X^{17} + 3 X^{16} + 174 X^{15} - 357 X^{14} - 1083
       X^{13} + 3463 X^{12} \\ + 2001 X^{11} - 13218 X^{10} + 3150 X^9
       + 22479 X^8 - 14883 X^7 \\ - 15063 X^6 + 16155 X^5 + 741 X^4 -
       5073 X^3 + 1557 X^2 - 37
     \end{multline*}\ \vskip-12pt
   }
   & $(6,3,1)(6,3,1)$
   & $5$
   & $2$
   & $U$
    \\
\hline
5  \phv{11}{5}
   & $3$
   & $44$
   & $I$
   & $9\cO_k$
   & $C_{12} \times C_2$
   & $C_{12}$
   & \parbox{10cm}{
    \vskip-12pt \begin{multline*}
      X^{24} - 45 X^{22} + 801 X^{20} - 7460 X^{18} + 40758 X^{16} -
      137607 X^{14} \\ + 291465 X^{12} - 381516 X^{10} + 294180 X^8 \\
      - 122240 X^6 + 24288 X^4 - 2112 X^2 + 64
    \end{multline*}\ \vskip-12pt
    }
   & $(12,1,2)$
   & $3$
   & $38\,480$
   & $p$
    \\
\hline
6  \phv{11}{5}
   & $3$
   & $53$
   & $I$
   & $90\cO_k$
   & $C_6 \times C_3$
   & $C_3 \times C_3$
   & \parbox{10cm}{
    \vskip-12pt \begin{multline*}
      X^{18} - 9 X^{17} - 33 X^{16} + 462 X^{15} - 147 X^{14} - 7521
      X^{13} + 11377 X^{12} \\ + 47199 X^{11} - 92040 X^{10} - 144180
      X^9 + 293583 X^8 \\ + 245031 X^7 - 406925 X^6 - 245451 X^5 +
      209247 X^4 \\ + 118989 X^3 - 21813 X^2 - 14520 X - 1331
    \end{multline*}\ \vskip-12pt
    }
   & $(6,1,3)$
   & $5$
   & $1$
   & $U$
   \\
\hline
\end{tabular}
}
\end{sidewaystable}

\end{document}